\documentclass[a4paper]{article}
\usepackage[left=3.5cm,right=3.0cm,top=2.5cm,bottom=2.5cm]{geometry}
\usepackage{xcolor}
\usepackage{graphicx}
\usepackage{cite}
\usepackage{bm}
\usepackage{amsmath,amssymb,amscd,amsxtra,amsfonts}
\usepackage{lipsum}
\usepackage{amsfonts}
\usepackage{graphicx}
\usepackage{multirow}

\usepackage{array}
\usepackage{amssymb}
\usepackage{verbatim}

\usepackage{xcolor}

\usepackage{epstopdf}
\usepackage{amsmath,bm,stmaryrd}

\usepackage{amsthm}

\usepackage{mathtools}
\usepackage{indentfirst}

\usepackage{float}
\usepackage{diagbox}

\newcommand\keywordsname{Key words}
\newcommand\AMSname{AMS subject classifications}

\newenvironment{@abssec}[1]
{\if@twocolumn
\section*{#1}%
\else
\vspace{.05in}\footnotesize
\parindent .2in
{\upshape\bfseries #1. }\ignorespaces
\fi}

{\if@twocolumn\else\par\vspace{.1in}\fi}

\newenvironment{keywords}{\begin{@abssec}{\keywordsname}}{\end{@abssec}}

\providecommand{\Div}{\operatorname{div}}          % Divergence
\providecommand*{\Dist}[2]{\operatorname{dist}({#1};{#2})}   % distance
\providecommand*{\Dist}[2]{\Dist{#1}{#2}}
\newcommand{\Vb}{{\mathbf{b}}}

\newcommand{\Vv}{{\mathbf{v}}}

\newcommand{\Bf}{{\boldsymbol{f}}}

\newcommand{\Bn}{{\boldsymbol{n}}}

\newcommand{\Br}{{\boldsymbol{r}}}

\newcommand{\Bu}{{\boldsymbol{u}}}
\newcommand{\Bv}{{\boldsymbol{v}}}
\newcommand{\Bw}{{\boldsymbol{w}}}
\newcommand{\Bx}{{\boldsymbol{x}}}

\newcommand{\BB}{{\boldsymbol{B}}}

\newcommand{\BD}{{\boldsymbol{D}}}

\newcommand{\BF}{{\boldsymbol{F}}}

\newcommand{\BH}{{\boldsymbol{H}}}
\newcommand{\BI}{{\boldsymbol{I}}}
\newcommand{\BJ}{{\boldsymbol{J}}}
\newcommand{\BK}{{\boldsymbol{K}}}

\newcommand{\BM}{{\boldsymbol{M}}}

\newcommand{\BQ}{{\boldsymbol{Q}}}

\newcommand{\BU}{{\boldsymbol{U}}}
\newcommand{\BV}{{\boldsymbol{V}}}

\newcommand{\Cm}{\mathcal{M}}

\newcommand{\Ct}{\mathcal{T}}

\newcommand{\bbA}{\mathbb{A}}

\newcommand{\bbC}{\mathbb{C}}

\newcommand{\bbI}{\mathbb{I}}

\newcommand{\bbL}{\mathbb{L}}

\newcommand{\bbN}{\mathbb{N}}

\newcommand{\bbP}{\mathbb{P}}

\newcommand{\bbU}{\mathbb{U}}

\newcommand{\bbZ}{\mathbb{Z}}

\newcommand{\be}{\begin{eqnarray}}
\newcommand{\ee}{\end{eqnarray}}

\newcommand{\ben}{\begin{eqnarray*}}
\newcommand{\een}{\end{eqnarray*}}

\newtheorem{theorem}{Theorem}[section]%  meant for continuous numbers
\newtheorem{remark}{Remark}%
\newtheorem{algorithm}[theorem]{\bf Algorithm}

\numberwithin{equation}{section}

\title{Robust augmented Block Preconditioning for two-phase Incompressible Inductionless MHD}
\author{Bijie Song \thanks{School of Mathematics and Statistics, Henan University, Kaifeng 475004, China.}
\and
Lingxiao Li
\thanks{School of Mathematics and Statistics,
Henan University, Kaifeng 475004, China. (lxli@henu.edu.cn)}
\and
Jie Zhang \thanks{State Key Laboratory for Strength and Vibration of Mechanical Structures, School of Aerospace, Xi'an Jiaotong University, Xi'an 710049, China.}}
\begin{document}
\date{}
\maketitle

\begin{abstract}
  In this paper, we develop a robust augmented block preconditioner for three-dimensional, two-phase, incompressible inductionless MHD equations. We employ a high-order monolithic and charge-conservative arbitrary Lagrangian--Eulerian (ALE) finite element method as base scheme, for which the solution variables are current density, velocity, electric potential and pressure. Grad-div stabilization is added in both the current density equation and discrete momentum equation. Backward Euler method is used for the time discretization. At each time step, fully coupled solution strategy is adopted. This leads to large-scale, ill-conditioned, nonsymmetric, and indefinite dual saddle-point system. To solve this system efficiently, we construct an augmented block preconditioner within the constraint preconditioning framework. The ultimate preconditioner combines the augmented Lagrangian technique with approximate block factorization for  current density sub-system and fluid sub-system. The fully coupled four-field system is then solved using preconditioned FGMRES iteration. Numerical experiments confirm the high-order accuracy and discrete charge conservation of the ALE-FEM. Drop of heavy fluid in light fluid is taken as benchmark to demonstrate the robustness of the preconditioner with respect to the discontinuity of physical parameters and the mesh size.
\end{abstract}


\section{Introduction}\label{sec:introduction}
Magnetohydrodynamics (MHD) describes the interaction between electrically conducting fluids and electromagnetic fields and has been widely applied to liquid-metal systems, including liquid-metal blankets in magnetic confinement fusion, electromagnetic metallurgy, and liquid-metal transport \cite{Gerbeau2006,kho23,urgorrippcf}.
In many of these applications, conducting fluids contain bubbles or droplets, or possess free surfaces, giving rise to multiphase MHD flows with moving interfaces \cite{zhang14,zha18}. Depending on whether the magnetic field induced by the fluid motion is taken into account, MHD models can generally be classified as induction and inductionless models. For liquid-metal flows with a small magnetic Reynolds number, the induced magnetic field is negligible compared with the imposed magnetic field, leading to the inductionless MHD model \cite{Gerbeau2006,ni12,LiSIAM2019,urgorrippcf}.
In this paper, we mainly study the efficient and robust preconditioner for the finite element discretization of a three-dimensional, two-phase, incompressible inductionless MHD model with moving interface.

Numerous numerical methods have been developed for two-phase MHD models with moving interfaces. According to whether the computational mesh moves with the interface, these methods can be roughly divided into fixed Eulerian mesh methods and moving mesh methods. Fixed Eulerian mesh methods describe the interface evolution on a fixed computational mesh. Typical examples include the volume-of-fluid (VOF) method, the level-set method, and the phase-field method. These methods are flexible for large deformations and complex topological changes of the interface \cite{Abels2012,Ki2010,zha18}.
In addition, moving-interface problems can also be solved using unfitted finite element methods, such as CutFEM, on fixed background meshes \cite{Frachon2019}. The interface can cut through the background mesh, so the mesh does not need to conform to the interface. In contrast, moving mesh methods directly track material interfaces through the motion of the computational mesh. Typical approaches include Lagrangian methods and arbitrary Lagrangian--Eulerian (ALE) methods \cite{ger03,wu2014jcp,sun2021}. Moving-mesh methods have also been used for MHD problems. Han and Tang developed an adaptive moving-mesh method for two-dimensional ideal MHD using magnetic-potential formulation \cite{Han2007}. In a Lagrangian method, the computational mesh moves with the fluid particles and therefore naturally follows the material interface. However, severe mesh distortion may occur under large deformations.
ALE methods allow the mesh to move with a velocity different from the fluid velocity, which provides greater flexibility in mesh motion \cite{Anderson2018,Rieben2007}. ALE methods can be broadly divided into monolithic ALE methods and split ALE methods. In a monolithic ALE method, the ALE equations are solved in a single stage. In a split ALE method, the computation is typically separated into different stages, such as Lagrangian evolution, mesh optimization, and field remapping \cite{Anderson2018}. Recently, high-order curvilinear meshes have been used in high-order Lagrangian \cite{Cheng2008,dob12,nik22} and ALE methods \cite{Anderson2018} to improve the geometric approximation of moving interfaces and the spatial discretization accuracy. In this paper, we use a high-order monolithic ALE finite element method as the basic scheme, and then develop an efficient and robust augmented block preconditioning technique for the two-phase incompressible inductionless MHD equations.

This work uses the current density, fluid velocity, electric potential, and pressure as the primary unknowns and solves them simultaneously using a fully coupled strategy. It extends our recent work on the two-phase incompressible MHD model \cite{Li2026}. The high-order monolithic ALE finite element discretization gives rise to a large-scale, nonsymmetric, ill-conditioned, and indefinite double saddle-point system. Its efficient solution is challenging because the ALE convection term introduces nonsymmetry, while large contrasts in material parameters across the interface can cause severe ill-conditioning.

For single-phase MHD models, saddle-point linear systems arising from finite element discretizations and their preconditioning methods have been extensively studied. The theoretical foundations of the relevant methods mainly include block factorizations and Schur complement approximations
\cite{Benzi2005,Elman2014,Phillips2014,Phillips2016}, approximate commutator methods \cite{Elman2006}, augmented Lagrangian methods \cite{BenziOlshanskii2006,BenziWang2011,LiJCP2017},
operator preconditioning \cite{Mardal2011},
constraint preconditioning \cite{KGW2000,C2002,Li2022}, and multigrid methods \cite{Adler2016}.
For MHD problems, Shadid et al. developed a fully implicit, fully coupled Newton--Krylov--AMG solution framework for incompressible resistive MHD \cite{Shadid2010,Shadid2016}.
Phillips and collaborators \cite{Phillips2014,Phillips2016} developed robust block preconditioners for incompressible resistive MHD equations based on least-squares commutators.
Cyr, Badia, Wathen, and their collaborators exploited the physical block structure of the systems and Schur complement approximations to develop block preconditioners for reduced resistive MHD \cite{Cyr2013},
inductionless MHD \cite{Badia2014},
and mixed finite element discretizations of full MHD \cite{Wathen2017,Wathen2020}, respectively.
Laakmann et al. further developed an augmented Lagrangian preconditioner applicable to high Reynolds and coupling numbers \cite{Laakmann2022}.
Recently, Ohm et al. \cite{Ohm2024} devised
scalable block preconditioning for low Mach number compressible MHD
for problems in magnetic confinement fusion.
By exploiting the block structure of the discrete systems and the properties of the underlying continuous operators, and by constructing readily invertible block operators, Schur-complement based approximations \cite{Phillips2016},
or multilevel solvers \cite{Adler2016,Muralikrishnan2023},
can effectively improve the efficiency and robustness of Krylov subspace methods for single-phase coupled fluid--electromagnetic systems.

However, efficient preconditioning is more challenging for multiphase flow problems. The density, viscosity, and other material parameters differ between fluids. They may also exhibit large jumps across interfaces. These coefficient variations affect the discrete operators and their Schur complements approximation.
Preconditioners for single-phase problems with constant coefficients may therefore lose robustness. For the two-phase incompressible Navier--Stokes equations, variable density and viscosity have received considerable attention. He and Neytcheva \cite{He2012} investigated augmented-Lagrangian-based preconditioning methods for the discrete Oseen problem with variable viscosity and demonstrated through theoretical analysis and numerical experiments that some preconditioning techniques widely used for constant-viscosity Oseen problems remain effective when the viscosity varies. Bootland et al. extended the approximate-commutator-based PCD and LSC preconditioners to two-phase incompressible Navier--Stokes systems with variable density and viscosity. By incorporating phase-dependent physical coefficients into the Schur complement approximations, they improved the performance of the preconditioners under mesh refinement and large jumps in density and viscosity and demonstrated their effectiveness using dam-break problem \cite{Bootland2018,Bootland2019}. These studies highlight the importance of accounting for jumps in material parameters, such as density and viscosity, across phase interfaces when designing preconditioners.

To our knowledge, efficient preconditioned iterative methods for fully coupled multiphase MHD remain less studied than those for single-phase MHD and two-phase NS problems. Motivated by augmented Lagrangian preconditioners for incompressible Navier--Stokes and Oseen problems \cite{BenziOlshanskii2006,BenziWang2011,He2018}, we introduce augmentation techniques into the design of our preconditioner. On the other hand, constraint preconditioning provides an important framework for indefinite linear systems with explicit constraint structures \cite{KGW2000,C2002}. The framework preserves the constraint blocks of the original system. It controls the spectral properties of the preconditioned system through a suitable approximation of the leading block on the constraint null space. Li et al. \cite{Li2022} applied this idea to the dual saddle-point system arising from the induction equation of resistive MHD using current density and magnetic potential
as solution variables. They combined augmentation techniques with approximate block factorizations to construct a robust block preconditioner.

Motivated by the above studies, we construct a robust block preconditioner for the four-field dual saddle-point system arising from a high-order ALE mixed finite element discretization of two-phase inductionless MHD. The preconditioner is developed within the constraint preconditioning framework and combines augmentation techniques with approximate block factorizations.
We further develop a preconditioned FGMRES solver using the proposed preconditioner.
The main contributions of this work are summarized as follows:
\begin{itemize}
   \item We construct a fully discrete high-order ALE finite element scheme for the 3D two-phase inductionless MHD equations. High-order curvilinear meshes are used to represent the interface \cite{Anderson2018,nik22}.
       The current density $\BJ$ is approximated by parametric $H(\mathrm{div})$-conforming BDM finite elements, which, together with the discrete charge-conservation equation, ensures discrete charge conservation, $\nabla\cdot\BJ_h=0$. The velocity and pressure are discretized by high-order parametric Taylor--Hood finite elements, and the electric potential is approximated by parametric discontinuous finite elements.
   \item For the above dual saddle-point system, we construct a block preconditioner within the constraint preconditioning framework by combining augmentation techniques with approximate block factorizations. For the current-density subsystem, we use a continuous augmented term in the variational formulation, and devise an approximation of the electric-potential Schur complement. For the fluid subsystem, a weighted pressure Schur complement approximation is constructed from the velocity mass matrix and grad-div terms in the
       momentum equation, incorporating the information of density, viscosity and the grad-div parameter. Extensive numerical experiments indicate that the block preconditioner is robust with respect to the discontinuity of main physical parameters and mesh refinement.
\end{itemize}
\vspace{2mm}

The structure of this article is organized as follows: Section~\ref{sec:Mathematical model} introduces the three-dimensional incompressible inductionless MHD model and derives its ALE strong formulation. In Section~\ref{sec:fully_discrete},
we discretize the model in space and time to obtain the fully discrete scheme, and present the detailed ALE algorithm. Section~\ref{sec:preconditioner} analyzes the dual saddle-point block structure of the fully discrete system and develops a block preconditioner within the constraint preconditioning framework by combining augmentation techniques and approximate block factorizations. In Section~\ref{sec:numerics}, a series of three-dimensional numerical experiments are presented to verify the spatial accuracy and charge-conservation property of the proposed method,
and in particular the parameters robustness and mesh refinement robustness of the preconditioner.
Finally, conclusions are given in Section~\ref{sec:conclusion}.

\section{Mathematical model and ALE strong formulation}\label{sec:Mathematical model}
%%%%%%%%%%%%%%%%%%%%%%%%%%%%%%%%%%%%%%%%%%%%%%%%%%%%%%%
\subsection{Two-phase incompressible  Inductionless MHD}

Let $\Omega(t)\subset\mathbb{R}^3$ be a bounded Lipschitz domain representing the current fluid region at
time $t \in [0,T_f]$. The domain is occupied by two immiscible incompressible fluids and is decomposed as
\[\Omega(t)=\Omega_+(t)\cup\Gamma(t)\cup\Omega_-(t),
   \qquad
   \Gamma(t)=\overline{\Omega_+(t)}\cap\overline{\Omega_-(t)},
   \qquad
   \Omega_{\pm}:=\Omega_+(t)\cup\Omega_-(t).\]
where $\Gamma(t)$ denotes the material interface separating the two phases. In the typical two-phase setting, we assume that the light liquid occupies the domain $\Omega_-(t)$, while the heavy fluid occupies the domain $\Omega_+(t)$.

We use the characteristic length $L_c$, velocity $u_c$, density $\rho_c$, dynamic viscosity $\nu_c$, electrical conductivity $\sigma_c$, and magnetic flux density $B_c$ to nondimensionalize the MHD equations. We assume that the density, viscosity, and electrical conductivity are piecewise constants in $\Omega(t)$; that is, they are constant within each phase, $\Omega_+(t)$ or $\Omega_-(t)$, but may take different values in the two phases,
\begin{equation}\label{sec2.1: eq-1_physical_parameters}
   \rho(\Bx) =
   \begin{cases}
      \rho_+, & \text{in } \Omega_+(t),\\
      \rho_-, & \text{in } \Omega_-(t),
   \end{cases}
   ~~
   \nu(\Bx) =
   \begin{cases}
      \nu_+, & \text{in } \Omega_+(t),\\
      \nu_-, & \text{in } \Omega_-(t),
   \end{cases}
   ~~
   \sigma(\Bx) =
   \begin{cases}
      \sigma_+, & \text{in } \Omega_+(t),\\
      \sigma_-, & \text{in } \Omega_-(t).
   \end{cases}
\end{equation}
The minimum and maximum values of these dimensionless physical parameters satisfy the following inequalities:
\begin{subequations}\label{sec2.1: eq0_physical_parameters}
\begin{align}
    0 < \rho_{\min} \leq \rho(\Bx) \leq \rho_{\max} < \infty,\\
    0 < \nu_{\min} \leq \nu(\Bx) \leq \nu_{\max} < \infty,\\
     0 < \sigma_{\min} \leq \sigma(\Bx) \leq \sigma_{\max} < \infty.
\end{align}
\end{subequations}

The applied magnetic field $\boldsymbol{B}$ is prescribed. The unknowns are the current density $\boldsymbol{J}$, fluid velocity $\boldsymbol{u}$, electric potential $\phi$ and pressure $p$ \cite{ni12,zha18,LiSIAM2019}. Then dimensionless two-phase incompressible inductionless MHD model equations
(see \cite{zha18,Li2026} for details)
\begin{subequations}
   \label{eq:model_1st}
   \begin{alignat}{2}
      {\sigma}^{-1}\boldsymbol{J}+\nabla\phi-\boldsymbol{u}\times\boldsymbol{B}
      &=\boldsymbol{0}, &\quad& \text{in } \Omega_{\pm}(t),
      \label{eq:model_1st:a}\\
      \rho\frac{d\boldsymbol{u}}{dt}-\nabla\cdot\boldsymbol{\tau}
      -\alpha\BJ\times\BB
      &=\rho\boldsymbol{f}_u, &\quad& \text{in } \Omega_{\pm}(t),
      \label{eq:model_1st:b}\\
      \nabla\cdot\boldsymbol{J}
      &=0, &\quad& \text{in } \Omega_{\pm}(t),
      \label{eq:model_1st:c}\\
      \nabla\cdot\boldsymbol{u}
      &=0, &\quad& \text{in } \Omega_{\pm}(t).
      \label{eq:model_1st:d}
   \end{alignat}
\end{subequations}
where $\boldsymbol{f}_u$ denotes the body force per unit mass, and $\frac{d}{dt}(\cdot)$ denotes the material derivative, defined for any sufficiently smooth variable $\psi(\cdot,t):\Omega(t)\rightarrow\mathbb{R}^3$ by
\begin{equation}\label{sec2.1: eq2_material_derivative}
   \frac{d\psi}{dt} :=  \frac{\partial\psi}{\partial t} + \boldsymbol{u}\cdot\nabla\psi.
\end{equation}

Besides, $\boldsymbol{\tau}$ denotes the Cauchy stress tensor
\[
   \boldsymbol{\tau}:=-p\bbI+\frac{2\nu}{Re}\boldsymbol{\varepsilon}(\boldsymbol{u}),
   \quad
   \boldsymbol{\varepsilon}(\boldsymbol{u}):=\frac{1}{2}\left(\nabla\boldsymbol{u}+\nabla\boldsymbol{u}^{T}\right),
   \quad
   \text{in }\Omega_{\pm}(t).
\]
where $\bbI\in\mathbb{R}^{3\times 3}$ denotes the identity matrix, $\boldsymbol{\varepsilon}(\boldsymbol{u})$ is the rate-of-deformation tensor, $Re = \frac{\rho_cu_cL_c}{\nu_c}$ is the Reynolds number
and $\alpha = \frac{\sigma_cB_c^{2}L_c}{\rho_cu_c}$ is the Stuart number.

We assume that $\Gamma(t)$ is a sufficiently smooth hypersurface. Let $\mathcal{V}$ denote the normal velocity of the interface and let $\Bn_{\Gamma}$ be the unit normal vector on $\Gamma(t)$ pointing from $\Omega_-(t)$ to $\Omega_+(t)$. Since $\Gamma(t)$ is a material interface, its normal velocity coincides with the normal component of the fluid velocity
\begin{equation*}
   \mathcal{V}=\boldsymbol{u}\cdot\Bn_{\Gamma} \qquad \text{on }\Gamma(t).
\end{equation*}
The following interface conditions are imposed on $\Gamma(t)$:
\begin{align*}
   [\boldsymbol{u}]_-^+ &= \boldsymbol{0},\\
   [\boldsymbol{\tau}\cdot\Bn_{\Gamma}]_-^+ &= \boldsymbol{0},\\
   [\phi]_-^+ &= 0,\\
   [\boldsymbol{J}\cdot\Bn_{\Gamma}]_-^+ &= 0.
\end{align*}
Here, $[\cdot]_-^+$ denotes the jump across $\Gamma(t)$, defined by
\[
   [q]_-^+ := q_+ - q_-,
\]
where $q_{\pm}$ denote the traces of $q$ on $\Gamma(t)$ taken from $\Omega_{\pm}(t)$, respectively.
For the concise presentation, we assume the following no-slip boundary condition for velocity and insulating boundary condition for current density,
\begin{equation*}
\boldsymbol{u} = \boldsymbol{0}, \quad \BJ\cdot \Bn_{\partial\Omega} = 0 \qquad \text{on}~~\partial\Omega(t),
\end{equation*}
where $\Bn_{\partial\Omega}$ denotes the outward unit normal vector on $\partial\Omega$. The initial velocity is prescribed as
\begin{equation*}
 \boldsymbol{u}(\cdot,0) = \boldsymbol{u}_0\qquad \text{on}~~\Omega(0).
\end{equation*}

%%%%%%%%%%%%%%%%%%%%%%%%%%%%%%%%%%%%%%%%%%%%%%%%%%%
\subsection{Strong formulation in the ALE framework}
Let $\Omega_0 := \Omega(0)$ be the ALE reference domain associated with $\Omega(t)$. For $t\in[0,T_f]$, the ALE mapping is defined by \cite{ger03,sun2021}
\begin{equation} \label{sec2.2: eq1_ALE_mapping}
   \Cm(\cdot,t):\Omega_0\rightarrow\Omega(t),
   \qquad
   {\Bx_0}\mapsto\Bx=\Cm({\Bx_0},t).
\end{equation}
where ${\Bx_0} \in \Omega_0$ denotes the ALE reference coordinate. Then ALE mesh velocity is defined by
\begin{equation}\label{sec2.2: eq2_ALE_velocity}
   \boldsymbol{w}(\boldsymbol{x},t) := \frac{\partial \Cm(\widehat{\Bx},t)}{\partial t},
   \qquad
   \boldsymbol{x} = \Cm(\widehat{\Bx},t)\in\Omega(t).
\end{equation}
We remark that the ALE mapping is not unique. In this work, the normal component of the ALE mesh velocity coincides with that of the fluid velocity on the material interface and the physical boundary
\begin{equation}\label{sec2.2: eq3_bdy}
   (\boldsymbol{w}-\boldsymbol{u})\cdot\Bn_{\Gamma} = 0 \quad \text{on}~~~\Gamma(t),
   \qquad
   (\boldsymbol{w}-\boldsymbol{u})\cdot\Bn_{\partial\Omega} = 0 \quad \text{on}~~~\partial\Omega(t).
\end{equation}
Following the standard ALE framework, for any sufficiently smooth field $\psi$ defined on the evolving domain, we define the ALE time derivative with respect to the ALE velocity $\Bw$ by
\begin{equation}\label{sec2.2: eq4_ALE_derivative}
   \frac{D\psi}{D t} = \frac{\partial\psi}{\partial t}  +  \boldsymbol{w}\cdot\nabla\psi.
\end{equation}
{Consequently, by \eqref{sec2.2: eq4_ALE_derivative}, the material derivative defined in \eqref{sec2.1: eq2_material_derivative} can be written in the ALE framework as}
\begin{equation} \label{sec2.2: eq5:ALE_material_derivative}
   \frac{d\psi}{dt} = \frac{D\psi}{D t} + \left(\boldsymbol{u}-\boldsymbol{w}\right) \cdot \nabla \psi.
\end{equation}
Using \eqref{sec2.2: eq5:ALE_material_derivative}, the strong form of system \eqref{eq:model_1st} in the ALE framework is given by
\begin{subequations}
   \label{eq:model_ALE}
   \begin{align}
      \sigma^{-1}\boldsymbol{J}+\nabla\phi-\boldsymbol{u}\times\boldsymbol{B}
      &=\boldsymbol{0}, &\quad& \text{in }\Omega_{\pm}(t),
      \label{eq:model_ALE:a}\\
      \rho\left[\frac{D\boldsymbol{u}}{Dt}
      +(\boldsymbol{u}-\boldsymbol{w})\cdot\nabla\boldsymbol{u}\right]
      -\nabla\cdot\boldsymbol{\tau}-\alpha\boldsymbol{J}\times\boldsymbol{B}
      &=\rho\boldsymbol{f}_u, &\quad& \text{in }\Omega_{\pm}(t),
      \label{eq:model_ALE:b}\\
      \nabla\cdot\boldsymbol{J}
      &=0, &\quad& \text{in }\Omega_{\pm}(t),
      \label{eq:model_ALE:c}\\
      \nabla\cdot\boldsymbol{u}
      &=0, &\quad& \text{in }\Omega_{\pm}(t).
      \label{eq:model_ALE:d}
   \end{align}
\end{subequations}

%%%%%%%%%%%%%%%%%%%%%%%%%%%%%%%%%%%%%%%%%%%%%%%%%%
\section{A high-order charge-conservative ALE mixed finite element method}\label{sec:fully_discrete}
In this section, we will give our base scheme, a high-order charge-conservative ALE-FEM, for which we will devise an augmented block preconditioner and finally develop a preconditioned FGMRES iterative solver for the two-phase
inductionless MHD.
The present scheme is an extension of our previous high-order charge-conservative Lagrangian FEM in \cite{Li2026}.

\subsection{Temporal discretization and high-order parametric mixed finite element spaces}
We uniformly partition the time interval $[0,T_f]$ as $0=t^0<t^1<\cdots<t^N=T$, where $t^n=n\Delta t$ and $\Delta t=T_f/N$.
Let $\Omega_h^0$ be a degree-$(k+1)$ isoparametric approximation \cite{BrennerScott2008}
of the initial physical domain $\Omega_0$, where $k \geq 1$. The $\Omega_h^0$ can be denoted by
\[
   \Omega_h^0=\Omega_{h,+}^0\cup\Gamma_h^0\cup\Omega_{h,-}^0.
\]
We denote by $\mathcal{T}_h^0$ a shape-regular and quasi-uniform interface-fitted curved tetrahedral mesh of $\Omega_h^0$ consisting of $N_e$ tetrahedral elements, such that $\overline{\Omega_h^0} = \bigcup_{K^0\in\mathcal{T}_h^0} \overline{K^0}$, where $K^0$ denotes a generic element of $\mathcal{T}_h^0$.
We define the mesh size $h$ as the maximum diameter of the initial straight-sided elements mesh $\widetilde{\mathcal{T}}_h^0$ by
\begin{equation}\label{sec3.1: eq2_mesh_size}
   h := \max_{\widetilde{K}^0\in\widetilde{\mathcal{T}}_h^0} h_{\widetilde{K}^0},\qquad
   h_{\widetilde{K}^0} :=\operatorname{diam}(\widetilde{K}^0),
\end{equation}
where $\widetilde{\mathcal{T}}_h^0$ denotes the underlying straight-sided partition for $\Omega_0$,
and each $\widetilde K^0\in\widetilde{\mathcal{T}}_h^0$ is the affine straight-sided element corresponding to $K^0\in\mathcal{T}_h^0$. Let $\boldsymbol{F}_{K^0}:\widehat{K}\rightarrow K^0$ denote the degree-$k$ isoparametric mapping from the reference tetrahedron $\widehat{K}$ onto $K^0$. We denote by $D\boldsymbol{F}_{K^0}$ its Jacobian matrix and define
\[
   \mathcal{J}_{K^0} := \det\left(D\boldsymbol{F}_{K^0}\right).
\]

The parametric finite element spaces on the initial mesh $\mathcal{T}_h^0$ (maybe curved) are defined by \cite{Li2026steady}
\begin{align}
   \BU_h^0 &:=
   \left\{\boldsymbol{\varphi}_h\in\boldsymbol{H}_0\bigl(\operatorname{div};\Omega_h^0\bigr):
   \left.\boldsymbol{\varphi}_h\right|_{K^0}\circ\boldsymbol{F}_{K^0}
   =\frac{D\boldsymbol{F}_{K^0}}{\mathcal{J}_{K^0}}\widehat{\boldsymbol{\varphi}},
   \forall~\widehat{\boldsymbol{\varphi}}\in\mathbf{BDM}_k(\widehat K),
   K^0\in\mathcal{T}_h^0\right\},
   \label{eq:initial_current_density_space}\\
   \widetilde{\BV}_h^0 &:=
   \left\{\boldsymbol{v}_h\in\boldsymbol{H}^1(\Omega_h^0):
   \left.\boldsymbol{v}_h\right|_{K^0}\circ\boldsymbol{F}_{K^0}
   =\widehat{\boldsymbol{v}},
   \quad \forall~\widehat{\boldsymbol{v}}\in\boldsymbol{P}_{k+1}(\widehat K),
   \quad K^0\in\mathcal{T}_h^0\right\},\\
   \BV_h^0 &:= \widetilde{\BV}_h^0\cap\BH_0^1(\Omega_h^0), \label{eq:initial_velocity_space}\\
   S_h^0 &:=
   \left\{r_h\in L^2(\Omega_h^0):
   \left.r_h\right|_{K^0}\circ\boldsymbol{F}_{K^0}
   =\widehat{r},\quad \forall~\widehat r\in P_{k-1}(\widehat K),
   \quad K^0\in\mathcal{T}_h^0\right\},
   \label{eq:initial_potential_space}\\
   Q_h^0 &:=
   \left\{q_h\in H^1(\Omega_h^0):
   \left.q_h\right|_{K^0}\circ\boldsymbol{F}_{K^0}
   =\widehat q,
   \quad \forall~\widehat q\in P_k(\widehat K),
   \quad K^0\in\mathcal{T}_h^0\right\},
   \label{eq:initial_pressure_space}
\end{align}
Here, $P_k(K^0)$ denotes the space of scalar polynomials of maximum total degree $k$ on $K^0$,
$\boldsymbol{P}_k(\hat{K}) := [P_k(\hat{K})]^3$. $\mathbf{BDM}_k(\widehat K)$ denotes the reference Brezzi--Douglas--Marini space of order $k$ \cite{boffibook}.

At each time level $t^n$, let $\BM_h^n \in \widetilde{\BV}_h^0$ (which will be decided in the following)
be the fully discrete ALE mapping which maps the initial configuration to the current computational domain.
Using $\BM_h^n$, for $n \geq 1$, we define the current curved mesh and domain by
\[\Ct_h^n := \bigl\{ K^n = \BM_h^n(K_0): K_0\in\Ct_{h}^0\bigr\},\quad \Omega_h^n:=\cup_{K^n\in \Ct_h^n} \overline{K}^n.\]
Denote the coordinates in $\Omega_h^n$ by $\Bx_h^n$. The ALE mapping in \eqref{sec2.2: eq1_ALE_mapping} is discretized at $t^n$ as
\begin{equation}\label{sec3.1: eq3_discretized_ALE_map}
    \BM_h^n:\Omega_h^0\rightarrow\Omega_h^n,\qquad
   \Bx_h^n=\BM_h^n(\Bx_h^0).
\end{equation}
The discrete interface $\Gamma_h^n :=\BM_h^n(\Gamma_h^0)$, is fitted to $\mathcal{T}_h^n$ and consists of faces of the curved tetrahedral elements. It separates  $\Omega_h^n$ into the two phase sub-domains.
\[
   \Omega_h^n=\Omega_{h,+}^n\cup\Gamma_h^n\cup\Omega_{h,-}^n.
\]
The ALE mesh motion is constructed so that this interface-fitted structure is preserved at each time level. Set
\[D\BM_h^n:=\nabla_{\Bx_h^0}\BM_h^n,\quad \mathcal{J}_{\BM_h^n}:=\det(D\BM_h^n).\]
Besides, we define $\boldsymbol{F}_{K^n}:=\BM_h^n\circ\BF_{K^0}:\widehat K\to K^n$ as the orientation-preserving element mapping. Let $D\boldsymbol{F}_{K^n}:=\nabla_{\widehat{\boldsymbol{x}}}\boldsymbol{F}_{K^n}$ and $\mathcal{J}_{K^n}:=\det(D\boldsymbol{F}_{K^n})$ denote the Jacobian matrix and the Jacobian determinant of $\boldsymbol{F}_{K^n}$, respectively.

Next, we introduce the definitions of the finite element spaces on the ALE mesh. At each time level $t^n$, let $\BU_h^n$, $\BV_h^n$, $S_h^n$, and $Q_h^n$ denote the finite element spaces on moving curved mesh $\Ct_h^n$ obtained from the corresponding initial spaces by the Piola and composition mappings, which are defined by \cite{Li2026}
\begin{align}
	\BU_h^n	&:=	\left\{	\boldsymbol{\varphi}_h^n	\in	\boldsymbol{H}_0
	\bigl(\operatorname{div};\Omega_h^n\bigr)
	:	\boldsymbol{\varphi}_h\circ\BM_h^n	=	\frac{D\BM_h^n}{\mathcal{J}_{\BM_h^n}}
	\boldsymbol{\varphi}_{h}^0,	\quad	\boldsymbol{\varphi}_{h}^0 \in\BU_{h}^0	\right\},\\
	\widetilde{\BV}_h^n	&:=	\left\{	\boldsymbol{v}_h^n
	\in	\boldsymbol{H}^1(\Omega_h^n)	:	\boldsymbol{v}_h^n\circ\BM_h^n
	=	\boldsymbol{v}_{h}^0,	\quad	\boldsymbol{v}_{h}^0\in\widetilde{\BV}_{h}^0
	\right\},\\
    \BV_h^n &:= \widetilde{\BV}_h^n \cap \BH_0^1(\Omega_h^n),\\
	S_h^n	&:=	\left\{	r_h^n\in L^2(\Omega_h^n)	:	r_h^n\circ\BM_h^n	=
	r_{h}^0,\quad r_{h}^0 \in S_{h}^0	\right\},\\
	{Q}_h^n	&:=	\left\{	q_h^n\in H^1(\Omega_h^n)	:	q_h^n\circ\BM_h^n	=
	q_{h}^0,\quad	q_{h}^0\in {Q}_{h}^0	\right\}.
\end{align}
The discrete counterparts of the material parameters in \eqref{sec2.1: eq-1_physical_parameters}
are defined as follows
\begin{align*}
   \rho_h^n &:= \rho_+\chi_{\Omega_{h,+}^n}+\rho_-\chi_{\Omega_{h,-}^n},\\
   \nu_h^n &:= \nu_+\chi_{\Omega_{h,+}^n}+\nu_-\chi_{\Omega_{h,-}^n},\\
   \sigma_h^n &:= \sigma_+\chi_{\Omega_{h,+}^n}+\sigma_-\chi_{\Omega_{h,-}^n},
\end{align*}
where $\chi_D$ denotes the characteristic function of a set $D\subseteq\Omega_h^n$.
{It follows from \eqref{sec2.1: eq0_physical_parameters} that the discrete physical parameters satisfy
\begin{subequations}\label{sec3.1: eq4_physical_parameters}
\begin{align}
    0 < \rho_{\min} \leq \rho_h^n \leq \rho_{\max} < \infty,\\
    0 < \nu_{\min} \leq \nu_h^n \leq \nu_{\max} < \infty,\\
     0 < \sigma_{\min} \leq \sigma_h^n \leq \sigma_{\max} < \infty.
\end{align}
\end{subequations}}
%%%%%%%%%%%%%%%%%%%%%%%%%%%%%%%%%%%%%%%%%%%%%%%%
\subsection{ALE mesh optimization} \label{sec:3.2}
%%%%%%%%%%%%%%%%%%%%%%%%%%%%%%%%%%%%%%%%%%%%%%%%%%
Denote by $\Upsilon_h^n = \Gamma_h^n \cup \partial\Omega_h^n$ the discrete curved surfaces consisting of the interface and the domain boundary.
For later use, we introduce the following notation for the $L^2$ inner products. Denote the $L^2$ inner product over the bulk domain $\Omega_{h,\pm}^n$ by
\[(u,v)_{\Omega_{h}^n} :=  \int_{\Omega_{h,+}^n}  u\cdot v\,d\boldsymbol{x}
   +  \int_{\Omega_{h,-}^n}    u\cdot v\,d\boldsymbol{x},\]
Firstly for $n\geq 1$, we introduce the homogeneous subspace
\[\mathcal{W}_{h,0}^{n-1}:=\left\{\boldsymbol{v}_h\in\widetilde{\boldsymbol{V}}_h^{n-1}:
   \boldsymbol{v}_h = \mathbf{0} \ \text{on }\Upsilon_h^{n-1}\right\}.\]
We now describe the construction of the ALE mesh from $t^{n-1}$ to $t^{n}$. Assume that the mesh $\mathcal{T}_h^{n-1}$ and the discrete physical fields at $t^{n-1}$ are known.
We firstly construct a predicted mesh position at $t^{n}$ using the fluid velocity at $t^{n-1}$
using Forward Euler's method
\[\bar{\Bx}_h^{n} = \boldsymbol{x}_h^{n-1}+\Delta t\,\boldsymbol{u}_h^{n-1}.\]
Then, for each $n\geq 1$, by the limited-harmonic mesh smoothing
\cite{Anderson2018},
we seek the grid functions $\boldsymbol{x}_h^{n}\in\widetilde{\BV}_h^{n-1}$ such that,
\begin{subequations}
   \begin{align}
   \omega\left(\nabla\boldsymbol{x}_h^{n},\nabla\boldsymbol{\Bv}_h\right)_{\Omega_{h}^{n-1}} +  \left(\Bx_h^n,\Bv_h\right)_{\Omega_{h}^{n-1}}
      &=\left(\bar{\boldsymbol{x}}_h^{\,n},\Bv_h\right)_{\Omega_{h}^{n-1}},
      \quad \forall\,\Bv_h\in\mathcal{W}_{h,0}^{n-1},\\
      \boldsymbol{x}_h^{n} &=\bar{\Bx}_h^{n}, \quad \text{on }\Upsilon_h^{n-1}.
   \end{align}
   \label{eq:harmonic_mesh_extension}
\end{subequations}
where $\omega>0$ is mesh smoothing parameter.
In the numerical tests, we always set $\omega = 10^{-4}$ in \eqref{eq:harmonic_mesh_extension}.
Once the updated mesh position has been obtained, the discrete mesh velocity is computed by
\begin{equation}\label{eq:discrete_mesh_velocity}
\Bw_h^{n-1}:=\frac{\Bx_h^n - \Bx_h^{n-1}}{\Delta t},\quad \mathrm{in}~~\Omega_h^{n-1}.
\end{equation}
For $n\geq1$, the ALE mapping between two consecutive meshes is defined by
\begin{equation}
   \BM_h^{n-1,n}: \Omega_h^{n-1}\rightarrow\Omega_h^n. \label{eq:consecutive_ALE_mapping}
\end{equation}
Once $\boldsymbol{w}_h^{n-1}$ has been obtained, then the discrete ALE map from the known geometry at $t^{n-1}$ to the updated geometry at $t^{n}$ will be
\begin{equation}
   \BM_h^{n-1,n}:= \operatorname{Id} + \Delta t\,\boldsymbol{w}_h^{n-1}.
   \label{eq:discrete_ALE_map}
\end{equation}
where $\operatorname{Id}$ is the identity mapping.
We denote the inverse of $\BM_{h}^{n-1,n}$ by $\BM_h^{n,n-1}$.
The discrete domain, interface, and phase subdomains then satisfy
\begin{equation}\label{eq:discrete_geometry_update}
   \Omega_h^{n}=\BM_h^{n-1,n}\left(\Omega_h^{n-1}\right),~~
   \Gamma_h^{n}=\BM_h^{n-1,n}\left(\Gamma_h^{n-1}\right),~~
   \Omega_{h,\pm}^{n}=\BM_h^{n-1,n}\left(\Omega_{h,\pm}^{n-1}\right).
\end{equation}
\begin{remark}
In fact, the mesh optimization strategy in this paper based on the limited-harmonic mesh smoothing is not optimal.
The mesh motion and optimization strategy is mainly intended to ensure the method is ALE.
In this work, we mainly focus on efficient and robust preconditioning techniques for the resulting
algebraic system.
More efficient mesh optimization techniques will be investigated in future work.
\end{remark}
%%%%%%%%%%%%%%%%%%%%%%%%%%%%%%%%%%%%%%%%%%%%%%%%%%
\subsection{Fully Discrete Scheme}
%%%%%%%%%%%%%%%%%%%%%%%%%%%%%%%%%%%%%%%%%%%%%%%%%%
{We now present the fully discrete ALE finite element scheme.} For $n=1,\ldots,N$, suppose that the discrete curved mesh $\Ct_h^n$ and the discrete fields $(\boldsymbol{u}_h^{n-1}, \Bw_h^{n-1})$ at $t^{n-1}$ are known.
On the updated mesh $\Ct_h^{n}$, using Backward Euler's method for the physical variables,
we seek $\left(\boldsymbol{J}_h^{n},\boldsymbol{u}_h^{n}, \phi_h^{n},p_h^{n}\right)\in\BU_h^{n}\times\BV_h^{n}\times S_h^{n}\times Q_h^{n}$ such that, for all $\left(\boldsymbol{\varphi}_h,\boldsymbol{v}_h,r_h,q_h\right)\in\BU_h^{n}\times\BV_h^{n}\times S_h^{n}\times Q_h^{n}$,
the following variational formulation holds:
\begin{subequations}
   \label{eq:fully_discrete_scheme}
   \begin{align}
   \left(\frac{1}{\sigma_h^{n}}\boldsymbol{J}_h^{n},
   \boldsymbol{\varphi}_h\right)_{\Omega_{h}^{n}} &+
   \gamma\left(\nabla\cdot\BJ_h^n,\nabla\cdot\boldsymbol{\varphi}_h\right)_{\Omega_{h}^{n}}
   -\left(\phi_h^{n},\nabla\cdot\boldsymbol{\varphi}_h\right)_{\Omega_{h}^{n}}\notag\\
   &-\left(\boldsymbol{u}_h^{n}\times\boldsymbol{B},
   \boldsymbol{\varphi}_h\right)_{\Omega_{h}^{n}}=0,\label{eq:fully_discrete_current}\\
   \frac{\left(\rho_h^{n}\Bu_h^{n},\Bv_h\right)_{\Omega_{h}^n}}{\Delta t}
   +&\left(\rho_h^n\left(\breve{\Bu}_h^{n-1}-\breve{\Bw}_h^{n-1}\right)\cdot\nabla\Bu_h^n,\Bv_h\right)_{\Omega_{h}^n}
   +\frac{1}{Re}\left(2\nu_h^{n}\boldsymbol{\varepsilon}(\Bu_h^{n}),\boldsymbol{\varepsilon}(\Bv_h)\right)_{\Omega_{h}^n}
   \notag\\
   &+\beta\left(\nabla \cdot\Bu_h^n,\nabla \cdot\Bv_h\right)_{\Omega_{h}^n}
   -\left(p_h^{n},\nabla\cdot\Bv_h\right)_{\Omega_{h}^n}
   -\alpha\left(\boldsymbol{J}_h^{n}\times\BB,\boldsymbol{v}_h\right)_{\Omega_{h}^n}\notag\\
   & =\frac{1}{\Delta t}\left(\rho_h^n\breve{\Bu}_h^{n-1},\Bv_h\right)_{\Omega_h^n}
   +\left(\rho_h^{n}\Bf_u^{n},\Bv_h\right)_{\Omega_{h}^{n}},\label{eq:fully_discrete_momentum}\\
   &-\left(\nabla\cdot\boldsymbol{J}_h^{n},r_h\right)_{\Omega_{h}^n}=0,\label{eq:fully_discrete_charge}\\
   &-\left(\nabla\cdot\boldsymbol{u}_h^{n},q_h\right)_{\Omega_{h}^n}=0.\label{eq:fully_discrete_incompressibility}
   \end{align}
\end{subequations}
where $0 < \beta = O(1)$ is the classical grad-div stabilization parameter \cite{or04,FAJN2017} and is one core ingredient in the preconditioning \cite{BenziOlshanskii2006,BenziWang2011,LiJCP2017}. Moreover, we added a similar augmented term
$\gamma\left(\nabla\cdot\BJ_h^n,\nabla\cdot\boldsymbol{\varphi}_h\right)_{\Omega_{h}^{n}}$ with $\gamma > 0$ in the
discrete variational equation of current density. Because $\Div\BJ_h=0$, this term will not change the discrete solutions.
However, $\beta\left(\nabla \cdot\Bu_h^n,\nabla \cdot\Bv_h\right)_{\Omega_{h}^n}$ will indeed change the discrete solutions.
So followed by the suggestions in \cite{or04,FAJN2017}, we restrict $\beta = O(1)$. The convective velocity is defined by
\begin{equation}
\breve{\Bu}_h^{n-1}:=\Bu_h^{n-1}\circ\BM_h^{n,n-1},\quad \breve{\Bw}_h^{n-1}:=\Bw_h^{n-1}\circ\BM_h^{n,n-1}.
\end{equation}
If $\Bw_h^{n-1} = \Bu_h^{n-1}$, then no mesh optimization is implemented and this scheme
reduces to the high-order charge-conservative Lagrangian FEM introduced in
\cite{Li2026}. Moreover, if $\Bw_h^{n-1} = \mathbf{0}$,
it reduces to the so-called semi-implicit Euler's method (with no material interface).
As described above, the corresponding mesh velocity $\boldsymbol{w}_h^{n-1}$ and the updated geometry $\Omega_h^{n}$
are firstly constructed. Therefore, when solving for the discrete fields at $t^{n}$, the relative convective velocity $\breve{\Bu}_h^{n-1}-\breve{\Bw}_h^{n-1}$ is already known, and hence the ALE convection term is linear with respect to the new velocity $\boldsymbol{u}_h^{n}$.
After solving \eqref{eq:fully_discrete_scheme}, the resulting velocity $\boldsymbol{u}_h^{n}$ is used to construct the predicted mesh velocity for the next time step,
and the above procedure in Subsection \ref{sec:3.2} is repeated until $t^n=T_f$.

%%%%%%%%%%%%%%%%%%%%%%%%%%%%%%%%%%%%%%%%%%%%%
%%%%%%%%%%%%%%%%%%%%%%%%%%%%%%%%%%%%%%%%%%%%%%%%%
\section{An Augmented Lagrangian Block Preconditioner}\label{sec:preconditioner}
The fully discrete ALE formulation \eqref{eq:fully_discrete_scheme} will result in a
coupled dual saddle-point linear system at each time step,
which is indefinite and generally ill-conditioned.
A direct solution of such a system becomes computationally expensive,
especially for large-scale three-dimensional problems.
Block preconditioning provides an effective framework for solving large-scale coupled saddle-point systems arising from multiphysics problems \cite{Elman2014}. Motivated by constraint preconditioning theory \cite{KGW2000,C2002,Li2022},
we construct a block preconditioner that exploits the constraint structure of the coupled system,
while augmentation and approximate block factorizations are further employed to obtain efficient approximations for the resulting saddle-point subproblems. Based on these ingredients, we develop an augmented block preconditioner and employ the preconditioned FGMRES method to solve the linear algebraic system arising from \eqref{eq:fully_discrete_scheme} at each time step.

Let $\{\boldsymbol{\varphi}_j^{n}\}_{j=1}^{N_{\BJ}}$, $\{\boldsymbol{v}_j^{n}\}_{j=1}^{N_{\Bu}}$, $\{r_j^{n}\}_{j=1}^{N_{\phi}}$, and $\{q_j^{n}\}_{j=1}^{N_{p}}$ be basis functions
of $\BU_h^{n}$, $\BV_h^{n}$, $S_h^{n}$, and $Q_h^{n}$, respectively.
The fully discrete system at the $t^{n}$ time level can be written in the algebraic form
\begin{equation}\label{sec4: eq1_algebraic_system}
\mathbb{A}^{n}\Vv^{n}=\Vb^{n},
\end{equation}
where
\begin{equation*}
   \Vv^{n} = \begin{pmatrix}\Vv_{\BJ}^{n}\\ \Vv_{\Bu}^{n}\\ \Vv_{\phi}^{n}\\ \Vv_{p}^{n}\end{pmatrix}.
\end{equation*}
contains the degrees of freedoms associated with $(\boldsymbol{J}_h^{n},\boldsymbol{u}_h^{n}, \phi_h^{n}, p_h^{n})$.
Then we get the coefficient matrix with the block structure
\begin{equation} \label{sec4: eq2_coefficient_matrix}
   \mathbb{A}^{n}
   =\begin{pmatrix}
      \widetilde{\BQ}_{\BJ}^{n} & (\BK^{n})^T & (G^{n})^T & 0\\
      -\alpha\BK^{n} & \BF^{n} & 0 & (B^{n})^T\\
      G^{n} & 0 & 0 & 0\\
      0 & B^{n} & 0 & 0
   \end{pmatrix}.
\end{equation}
where the correspond block matrices are defined by
{\begin{align*}
   \widetilde{\BQ}_{\BJ}^{n}[i,j]
   &:= \BQ_{\BJ}^{n}[i,j] + \gamma\left(\nabla\cdot\boldsymbol{\varphi}_j^{n},
   \nabla\cdot\boldsymbol{\varphi}_i^{n}\right)_{\Omega_{h}^{n}},\quad
   \BQ_{\BJ}^{n}[i,j]:=\left(\frac{1}{\sigma_h^{n}}\boldsymbol{\varphi}_j^{n},
   \boldsymbol{\varphi}_i^{n}\right)_{\Omega_{h}^n},\\
   \BK^{n}[i,j]
   &:=
   \left(\boldsymbol{\varphi}_j^{n}\times\boldsymbol{B},
   \boldsymbol{v}_i^{n}\right)_{\Omega_{h}^n},\\
   G^{n}[i,j]
   &:=
   -\left(\nabla\cdot\boldsymbol{\varphi}_j^{n},
   r_i^{n}\right)_{\Omega_{h}^n},\\
   B^{n}[i,j]
   &:=
   -\left(\nabla\cdot\boldsymbol{v}_j^{n},
   q_i^{n}\right)_{\Omega_{h}^n},\\
   \BF^{n}[i,j]
   &:=\frac{1}{\Delta t}\left(\rho_h^{n}\boldsymbol{v}_j^{n},\boldsymbol{v}_i^{n}\right)_{\Omega_{h}^n}
   +\left(\rho_h^n\left(\breve{\Bu}_h^{n-1}-\breve{\Bw}_h^{n-1}\right)\cdot\nabla\Bv_j^n,\Bv_i^n\right)_{\Omega_{h}^n}\notag\\
   &\quad+\frac{1}{Re}\left(2\nu_h^{n}\boldsymbol{\varepsilon}(\boldsymbol{v}_j^{n}),
   \boldsymbol{\varepsilon}(\boldsymbol{v}_i^{n})\right)_{\Omega_{h}^n}+\beta\left(\nabla \cdot\Bv_j^n,\nabla \cdot\Bv_i^n\right)_{\Omega_{h}^n}.
\end{align*}}
Based on the dual saddle-point system \eqref{sec4: eq2_coefficient_matrix}, we next derive the corresponding preconditioner within the constraint preconditioning framework, combined with augmented Lagrangian techniques. We fix an arbitrary time level $t^n$ and, for notational simplicity, suppress the superscript $n$ throughout this section unless otherwise stated.

%%%%%%%%%%%%%%%%%%%%%%%%%%%%%%%%%%%%%%%%%%%%%%%%%%%%%
\subsection{Constraint Preconditioning Framework}
In this section, we employ the constraint preconditioning framework introduced by Keller et al. for indefinite linear systems \cite{KGW2000}. The basic idea is to approximate the block associated with the primary physical variables while retaining the constraint blocks exactly. This framework was subsequently extended by Cao
to nonsymmetric indefinite systems \cite{C2002} and by Li et al. to single phase MHD equations in \cite{Li2022}.
In this work, we heuristically extend the constraint preconditioning framework
for two-phase inductionless MHD systems \eqref{eq:fully_discrete_scheme} and
further develop an augmented Lagrangian preconditioner.

Firstly, we write the coefficient matrix in \eqref{sec4: eq2_coefficient_matrix} in block form as follows:
\begin{equation}\label{sec4.1: eq1_block_A}
   \mathbb{A}=\begin{pmatrix}
   \bbZ & \bbN^T\\ \bbN & \boldsymbol{0}\end{pmatrix}.
\end{equation}
where the two block submatrices are defined by
\begin{equation}\label{sec4.1: eq2_Z_and_N}
   \bbZ:=\begin{pmatrix}\widetilde{\BQ}_{\BJ} & \BK^T\\ -\alpha \BK & \BF\end{pmatrix},
   \qquad
   \bbN:=\begin{pmatrix}G & 0\\ 0 & B\end{pmatrix}.
\end{equation}
Following the constraint preconditioning framework \cite{KGW2000,C2002,Li2022}, we retain the constraint coupling blocks $\bbN$ and $\bbN^T$ exactly and approximate $\bbZ$ by $\widetilde{\bbZ}$. If $\widetilde{\bbZ}$ provides a good
approximation to $\bbZ$ on $\ker(\bbN)$, then the resulting matrix $\bbP_e$ generally
is an effective preconditioner for $\mathbb{A}$.
The corresponding constraint preconditioner $\bbP_e$ is then defined by
\begin{equation}\label{sec4.1_eq3_constraint_preconditioner}
   \mathbb{P}_e=\begin{pmatrix}\widetilde{\bbZ} & \bbN^T\\ \bbN & \boldsymbol{0}\end{pmatrix}.
\end{equation}
We remark that on the kernel space of $\bbN$, the action of $\widetilde{\BQ}_{\BJ}$ on a vector with divergence-free
functions for current density block is equal to the action of ${\BQ}_{\BJ}$ on that vector. So define
\begin{equation}
   \widehat{\bbZ}:=\begin{pmatrix}{\BQ}_{\BJ} & \BK^T\\ -\alpha \BK & \BF\end{pmatrix}.
\end{equation}
and $\widehat{\bbZ}$ provides good approximation for $\bbZ$ on kernel space of $\bbN$.
To construct a good approximation $\widetilde{\bbZ}$ to $\bbZ$,
next we consider the exact block LU factorization of $\widehat{\bbZ}$
\begin{equation}\label{sec4.1_eq4_Z_LU}
   \widehat{\bbZ}=
   \begin{pmatrix}\BI & 0\\ -\alpha \BK\BQ_{\BJ}^{-1} & \BI\end{pmatrix}
   \begin{pmatrix}\BQ_{\BJ} & \BK^T\\ 0 & \BF+\alpha \BK\BQ_{\BJ}^{-1}\BK^T\end{pmatrix}.
\end{equation}
where $\BI$ is the identity matrix.

Based on the block  factorization in \eqref{sec4.1_eq4_Z_LU}, it suffices to construct a suitable approximation to the Schur complement $\BF+ \alpha \BK\BQ_{\BJ}^{-1}\BK^T$
to obtain an explicit form of $\widetilde{\bbZ}$. Since $\BQ_{\BJ}$ is the discrete representation of the linear operator $\frac{1}{\sigma_h^n}\operatorname{Id}$, we find that the inverse matrix $\BQ_{\BJ}^{-1}$ formally corresponds to the linear multiplication operator $\sigma_h^n\operatorname{Id}$. This motivates the following approximation of $\BF+ \alpha \BK\BQ_{\BJ}^{-1}\BK^T$:
\begin{align}\label{sec4.1: eq5_Fhat_approx}
   \widehat{\BF}[i,j]
   &:=\frac{1}{\Delta t}\left(\rho_h^{n}\boldsymbol{v}_j^{n},\boldsymbol{v}_i^{n}\right)_{\Omega_{h}^n}
   +\left(\rho_h^n\left(\breve{\Bu}_h^{n-1}-\breve{\Bw}_h^{n-1}\right)\cdot\nabla\Bv_j^n,\Bv_i^n\right)_{\Omega_{h}^n}\notag\\
   &\quad+\frac{1}{Re}\left(2\nu_h^{n}\boldsymbol{\varepsilon}(\Bv_j^{n}),\boldsymbol{\varepsilon}(\Bv_i^{n})\right)_{\Omega_{h}^n}+\beta\left(\nabla \cdot \Bv_j^n,\nabla \cdot \Bv_i^n\right)_{\Omega_{h}^n}\notag\\
   &+\alpha \left(\sigma_h^{n}\boldsymbol{B}\times\Bv_j^{n},\boldsymbol{B}\times\Bv_i^{n}\right)_{\Omega_{h}^{n}}.
\end{align}
Accordingly, we choose
\begin{equation} \label{sec4.1: eq6_Zwilde}
   \widetilde{\bbZ}=\begin{pmatrix}\widetilde{\BQ}_{\BJ} & \BK^T\\ 0 & \widehat{\BF}\end{pmatrix},
\end{equation}
where $\widehat{\BF}$ is an approximation to $\BF+\alpha\BK\BQ_{\BJ}^{-1}\BK^T$
defined in \eqref{sec4.1: eq5_Fhat_approx}.
Consequently, by \eqref{sec4.1: eq6_Zwilde} and \eqref{sec4.1_eq3_constraint_preconditioner},
we obtain the initial block preconditioner for the system $\mathbb{A}$ as
\begin{equation}\label{sec4.1: eq7_A1}
   \mathbb{A}_1 = \begin{pmatrix}
      \widetilde{\BQ}_{\BJ} & \BK^T & G^T & 0\\
      0 & \widehat{\BF} & 0 & B^T\\
      G & 0 & 0 & 0\\
      0 & B & 0 & 0
   \end{pmatrix}.
\end{equation}
Next, we reorder the unknowns as $(\boldsymbol{J}_h^{n},\phi_h^{n},\boldsymbol{u}_h^{n},p_h^{n})$ and rewrite $\mathbb{A}_1$ in the following block form:
\begin{equation}\label{sec4.1: eq8_A1_reorder}
   \mathbb{A}_1  = \begin{pmatrix}
      \mathbb{A}_{\BJ} & \bbC\\
      \boldsymbol{0} & \mathbb{A}_{\Bu}
   \end{pmatrix},
\end{equation}
where
\begin{equation} \label{sec4.1: eq9_AJ_Au}
   \mathbb{A}_{\BJ}:=\begin{pmatrix}\widetilde{\BQ}_{\BJ} & G^T\\ G & 0\end{pmatrix},\qquad \mathbb{A}_{\Bu}:=\begin{pmatrix}\widehat{\BF} & B^T\\ B & 0\end{pmatrix},\qquad
   \bbC:=\begin{pmatrix}\BK^T & 0\\ 0 & 0\end{pmatrix}.
\end{equation}
Consider the following residual equation associated with the preconditioner $\mathbb{A}_1$:
\begin{equation}\label{sec4.2: eq10_A1_residual}
   \mathbb{A}_1 \boldsymbol{\delta} = \boldsymbol{r},
\end{equation}
where $\boldsymbol{\delta} =(\delta_{\BJ},\delta_{\Bu},\delta_{\phi},\delta_p)^T$ denotes the correction vector,
and $\boldsymbol{r}=(r_{\BJ},r_{\Bu},r_\phi,r_p)^T$ denotes the residual vector. Then, using the block forms in \eqref{sec4.1: eq8_A1_reorder} and \eqref{sec4.1: eq9_AJ_Au}, we expand the system in \eqref{sec4.2: eq10_A1_residual} as follows:
\begin{subequations} \label{sec4.1: eq11_A1_residual_block}
   \begin{align}
      \mathbb{A}_{\BJ} \boldsymbol{\delta}_{1}+ \bbC \boldsymbol{\delta}_{2}
      &=\Br_1,
      \label{sec4.1: eq11a_A1_residual_block}\\
      \mathbb{A}_{\Bu} \boldsymbol{\delta}_2
      &=\Br_2,
      \label{sec4.1: eq11b_A1_residual_block}
   \end{align}
\end{subequations}
Therefore, solving the residual equation associated with the preconditioner $\mathbb{A}_1$ in \eqref{sec4.1: eq11_A1_residual_block} is equivalent to solving two constrained saddle-point subproblems: the current-density subsystem $\mathbb{A}_{\BJ}$ and the fluid subsystem $\mathbb{A}_{\Bu}$.
\vspace{2mm}

\begin{remark}
    In fact, although the preconditioner $\bbA_1$ is obtained by applying constraint preconditioning to $\bbA$ at the outermost level, the exact application of \(\bbA_1^{-1}\) in \eqref{sec4.2: eq10_A1_residual} remains computationally expensive. As shown in \eqref{sec4.1: eq9_AJ_Au} and \eqref{sec4.1: eq11_A1_residual_block}, this preconditioner gives rise to two indefinite saddle-point subsystems: the current-density saddle-point subsystem $\mathbb{A}_{\BJ}$ and the fluid saddle-point subsystem $\mathbb{A}_{\Bu}$. Therefore, suitable preconditioners must be constructed for these two saddle-point subsystems to solve them approximately.
    Based on this observation, the final preconditioner for the two-phase MHD model considered in this paper
    is constructed within the constraint preconditioning framework described above and then combines
    suitable preconditioning for the two saddle-point subsystems.
    Detailed derivation are provided in Sections~\ref{sec4.2:Current_Density_Subsystem} and~\ref{sec4.3: Fluid_Subsystem}.
\end{remark}
%%%%%%%%%%%%%%%%%%%%%%%%%%%%%%%%%%%%%%%%%%%%%%
\subsection{Preconditioning for the Current Density Subsystem}\label{sec4.2:Current_Density_Subsystem}
In this section, we consider the preconditioning for the current-density saddle-point subproblem $\bbA_{\BJ}$
in \eqref{sec4.1: eq9_AJ_Au}.
Our objective is to construct a robust and stable preconditioner for this saddle matrix.
To this end, we firstly apply precise LU decomposition to $\bbA_{\BJ}$ and then construct approximate Schur complement to the electric potential block.
Applying the LU decomposition of $\bbA_{\BJ}$ gives
\begin{align}
   \bbA_{\boldsymbol{J}} =
   \begin{pmatrix}
      \widetilde{\boldsymbol{Q}}_{\boldsymbol{J}} & {G}^T \\
      G & {0}
   \end{pmatrix} =
   \begin{pmatrix}
      \boldsymbol{I} & {0} \\
      G {\widetilde{\boldsymbol{Q}}}_{\boldsymbol{J}}^{-1} & \boldsymbol{I}
   \end{pmatrix}
   \begin{pmatrix}
      {\widetilde{\boldsymbol{Q}}}_{\boldsymbol{J}} & {G}^T \\
      {0} & -S_{\phi}
   \end{pmatrix} =:
   \bbL_{\boldsymbol{J}}\bbU_{\boldsymbol{J}}.
\end{align}
where $S_\phi$ is the exact Schur complement matrix of electrical potential $\phi$
\begin{equation}\label{sec4.2: eq3_Sphi}
   S_\phi := G {\widetilde{\BQ}}_{\BJ}^{-1}G^T.
\end{equation}

Following the framework in the monograph \cite{Elman2014}, it suffices to construct an approximate Schur complement that is computationally efficient and provides good approximation to $S_\phi$.
To do this, observing that when $\gamma = 0$, the precise Schur complement is $G {\BQ}_{\BJ}^{-1}G^T$.
Mass matrix is generally spectrally equivalent to its diagonal matrix.
Define $\BD_\BJ=\operatorname{diag}(\BQ_\BJ)$.
In this limiting case, we can use the approximate Schur complement of $S_{\phi}$ by
\begin{equation}
L_\phi:=G \BD_{\BJ}^{-1}G^T.
\end{equation}

When the electrical conductivity is constant or the variation is mild, and $\gamma > 0$,
the derivation in \cite{Arnold1997,Li2022,Vassilevski1996}
showed that
\begin{equation}
\frac{1}{\gamma} Q_\phi,\quad Q_\phi[i,j]: =\left(r_j^{n}, r_i^{n}\right)_{\Omega_{h}^n} .
\end{equation}
is good approximation.

Combining the two cases, we propose the following approximate Schur complement for $\phi$
\begin{equation}\label{sec4.2: eq6_hat_Sphi}
\widehat{S}_{\phi}^{-1} := L_{\phi}^{-1} + \gamma Q_{\phi}^{-1}.
\end{equation}
Finally we define the following preconditioner:
\begin{equation}\label{sec4.2: eq7_PJ}
   \mathbb{P}_{\BJ}=
   \begin{pmatrix}
      \widetilde{\boldsymbol{Q}}_{\BJ} & {G}^T \\ 0 & -\widehat{S}_{\phi}
   \end{pmatrix}.
\end{equation}

\begin{remark}
The information of discontinuity of $\sigma_h^n$ is embedded into the approximate Schur complement by $\BD_\BJ$.
$L_\phi$ is a scaled Poisson matrix, which can be efficiently solved by BoomerAMG preconditioner.
By our extensive numerical tests, letting $\gamma=0$ will obtain relative poor preconditioning than $\gamma \neq 0$.
Moreover, without $L_\phi$ and using only $\frac{1}{\gamma} Q_\phi$ as the approximation will get bad results
when $\sigma_h^n$
is discontinuous across the interface. Overall we find that the combination of $L_\phi$ and $\frac{1}{\gamma} Q_\phi$
is a robust choice.
\end{remark}

\subsection{Preconditioning for the Fluid Subsystem}\label{sec4.3: Fluid_Subsystem}
Next, we consider the construction of preconditioner for the fluid subsystem. By applying a block LU factorization to the saddle point system $\bbA_{\Bu}$, we obtain
\begin{equation}
   \bbA_{\Bu} = \begin{pmatrix}\BI & 0\\ B\widehat{\BF}^{-1} & \BI\end{pmatrix}
   \begin{pmatrix}\widehat{\BF} & B^T\\ 0 & -S_p\end{pmatrix},
\end{equation}
where
\begin{equation}
   S_p := B\widehat{\BF}^{-1}B^T
\end{equation}
is the pressure Schur complement. A direct treatment of this Schur complement requires repeated applications of $\widehat{\BF}^{-1}$, which is computationally expensive for large-scale three-dimensional problems. Therefore, in practical preconditioners, $S_p$ is usually replaced by an easily invertible approximation $\widehat S_p$ that preserves the essential spectral properties of $S_p$ \cite{Benzi2005,Elman2014,He2018}.

The continuous operator corresponding to $\widehat{\BF}$ is given by
\begin{align}
   \widehat{\mathcal F}\BU &=\frac{\rho_h^n}{\Delta t}\BU+\rho_h^n (\breve{\boldsymbol{u}}_h^{n-1}-\breve{\boldsymbol{w}}_h^{n-1})\cdot\nabla\BU - \frac{1}{Re}\nabla\cdot\bigl(2\nu_h^n\boldsymbol{\varepsilon}(\,\BU\,)\bigr)\notag\\
   &-\beta\nabla\operatorname{div}\BU - \alpha\sigma_h^n\boldsymbol{B}\times\left(\boldsymbol{B}\times\BU\right).
\end{align}
When the time derivative is neglected,
the operator $\widehat{\mathcal F}$ reduces to a
steady Oseen-type operator with some perturbation.
When stabilization term $\beta(\Div\Bu_h,\Div\Bv_h)$ is presented,
previous studies have shown that the pressure Schur complement can be
effectively approximated by a suitably weighted pressure
mass matrix \cite{BenziOlshanskii2006,BenziWang2011,He2012,LiJCP2017}.
Then, we introduce the weighted pressure mass matrix
\begin{equation}
M_p[i,j]:=\left(\frac{1}{\frac{\nu_h^n}{Re}+\beta}q_j^n,q_i^n\right)_{\Omega_h^n}.
\label{eq:weighted_pressure_mass}
\end{equation}
Therefore, in the steady-state case, $M_p$ can be used as an effective approximation
to the precise pressure Schur complement. Note that the non-uniform of $\nu_h^n$ is already accounted for
the computation in $M_p$.
Next, we consider the case with the time derivative. The matrix  $\widehat{\BF}$ can be written as
\begin{equation}
   \widehat{\BF}= \frac{1}{\Delta t}\BM_{\Bu} + \BF_0,
\end{equation}
where
\begin{align}
   \BM_{\Bu}[i,j]&:=\left(\rho_h^n\boldsymbol{v}_j^n,\boldsymbol{v}_i^n\right)_{\Omega_h^n},\\
   \BF_0[i,j] &:= \left(\rho_h^n\left(\breve{\Bu}_h^{n-1}-\breve{\Bw}_h^{n-1}\right)\cdot\nabla\Bv_j^n,\Bv_i^n\right)_{\Omega_{h}^n}
   +\frac{1}{Re}\left(2\nu_h^{n}\boldsymbol{\varepsilon}(\Bv_j^{n}),\boldsymbol{\varepsilon}(\Bv_i^{n})\right)_{\Omega_{h}^n}\\
   &+\beta\left(\nabla \cdot \Bv_j^n,\nabla \cdot \Bv_i^n\right)_{\Omega_{h}^n}
   +\alpha \left(\sigma_h^{n}\boldsymbol{B}\times\Bv_j^{n},\boldsymbol{B}\times\Bv_i^{n}\right)_{\Omega_{h}^{n}}.\notag
\end{align}
The precise pressure Schur complement can be rewritten by
\begin{equation}
   S_p=B\left(\frac{1}{\Delta t}\BM_{\Bu}+\BF_0\right)^{-1}B^T.
\end{equation}

When $\Delta t\rightarrow0$, transient mass term $\BM_\Bu$ dominates the velocity block. Formally, pressure Schur complement is expected to behave like $B\left(\frac{1}{\Delta t}\BM_{\Bu}\right)^{-1}B^T$.
Since the mass matrix $\BM_\Bu$ is spectrally equivalent to its diagonal approximation.
We then approximate this limiting operator by replacing $\BM_\Bu$ with its
diagonal approximation $\BD_\Bu=\operatorname{diag}(\BM_\Bu)$, and define
\begin{equation}
   L_p := B\BD_\Bu^{-1}B^T.
\end{equation}

Combining the above results, for sufficiently small time steps, $\frac{1}{\Delta t}L_p^{-1}$ provides an effective approximation to the ${S^{-1}_p}$. On the other hand, in the steady-state case without the time derivative, the weighted pressure mass matrix $M_p$ is an effective approximation
\cite{BenziOlshanskii2006,BenziWang2011,LiJCP2017}.
Motivated by these two limiting cases, we propose the following approximate Schur complement of $S_p$
\begin{equation}\label{sec4.3: eq_hatSp}
(\widehat{S}_p)^{-1}:= M_p^{-1} + \frac{1}{\Delta t}L_p^{-1}.
\end{equation}
We finally obtain the following preconditioner for matrix $\bbA_u$
\begin{equation}\label{sec4.3: eq_Pu}
\mathbb{P}_{\Bu} = \begin{pmatrix} \widehat{\BF} & B^T \\ 0 & -\widehat{S}_p \end{pmatrix}.
\end{equation}

\begin{remark}
In the devise of preconditioning above, the discontinuity of $\sigma_h^n$ is
taken into account in the computation of $\widehat{S}_{\phi}$ and $\widehat{\BF}$.
While the discontinuity of density $\rho_h^n$ and viscosity coefficient $\nu$ is also
embedded into the definition of disturbed transient reaction convection-diffusion matrix $\widehat{\BF}$.
Moreover, a core part is that
discontinuity of $\rho_h^n$ and $\nu_h^n$ should also be considered in the approximation of pressure schur complement,
which is the main philosophy for the pure Navier-Stokes equations \cite{He2012,Bootland2018,Bootland2019}.
\end{remark}
\vspace{2mm}

%%%%%%%%%%%%%%%%%%%%%%%%%%%%%%%%%%%%%%%%%%%
\subsection{Final Block Preconditioner}
%%%%%%%%%%%%%%%%%%%%%%%%%%%%%%%%%%%%%%%%%%%
We approximate the blocks $\bbA_{\BJ}$ and $\bbA_{\Bu}$ in \eqref{sec4.1: eq8_A1_reorder} by the preconditioners $\bbP_{\BJ}$ and $\bbP_{\Bu}$ defined in \eqref{sec4.1: eq9_AJ_Au} and \eqref{sec4.3: eq_Pu}, respectively.
Then, we obtain the final block upper triangular preconditioner for $\mathbb{A}$
in \eqref{sec4: eq2_coefficient_matrix}
\begin{equation} \label{sec4.4: eq_final_preconditioner}
   \mathbb{P}=
   \begin{pmatrix}
      \widetilde{\BQ}_{\BJ} & \BK^T & G^T & 0\\
      0 & \widehat{\BF} & 0 & B^T\\
      0 & 0 & -\widehat{S}_\phi & 0\\
      0 & 0 & 0 & -\widehat{S}_p
   \end{pmatrix},
\end{equation}
where $\widehat{S}_{\phi}^{-1}:= L_{\phi}^{-1} + \gamma Q_\phi^{-1}$ in \eqref{sec4.2: eq6_hat_Sphi}
and $\widehat{S}_p^{-1}:= M_p^{-1} + \frac{1}{\Delta t}L_p^{-1}$ in \eqref{sec4.3: eq_hatSp}.
For a given vector $\boldsymbol{r}=(r_\BJ,r_\Bu,r_\phi,r_p)^T$, the action of $\mathbb{P}^{-1}$ is obtained by solving
\begin{equation}\label{sec4.4: eq_Pdeta}
   \mathbb{P}\boldsymbol{\delta}=\boldsymbol{r}
\end{equation}
Owing to the block upper triangular structure of $\mathbb{P}$, the vector $\boldsymbol{\delta}$ can be computed by block substitution. Detailed solution procedure for \eqref{sec4.4: eq_Pdeta} is presented in Algorithm~\ref{sec4.4: algorithm}.

\begin{algorithm}\label{sec4.4: algorithm}
Solve $\mathbb{P}\boldsymbol{\delta}=\boldsymbol{r}$ in \eqref{sec4.4: eq_Pdeta} as follows:
\begin{enumerate}
   \item Compute $\delta_p=-\xi_p-\frac{1}{\Delta t}\eta_p$. We solve the following two subproblems

   \begin{itemize}
      \item Solve $ M_p \xi_p = r_p$ by PCG solver with diagonal preconditioning with relative tolerance is $10^{-4}$.
      The maximum iteration number is 30.
      \item Solve $L_p \eta_p = r_p$ by PCG solver with BoomerAMG preconditioning with relative tolerance is $10^{-4}$.
      The maximum iteration number is 30.
   \end{itemize}

   \item Solve $\widehat{\BF} \delta_\Bu=r_\Bu-B^T\delta_p$ using
   the GMRES solver with one level additive Schwarz preconditioner \cite{cai1999} with overlap 2.
   The relative residual tolerance is set to $10^{-4}$. The maximum iteration number is 30.

   \item Compute $\delta_\phi=-\xi_\phi- \gamma \eta_\phi$.
   We solve the following two subproblems

   \begin{itemize}
      \item Solve $L_\phi \xi_\phi = r_\phi$ by PCG solver with BoomerAMG preconditioning with
      relative tolerance is $10^{-4}$. The maximum iteration number is 30.
      \item Solve $Q_\phi \eta_\phi = r_\phi$ by PCG solver with diagonal preconditioning.
      The relative tolerance is $10^{-4}$ and maximum iteration number is 30.
   \end{itemize}

   \item Solve $\widetilde{\BQ}_\BJ \delta_\BJ=r_\BJ-\BK^T\delta_\Bu-G^T\delta_\phi$ using direct solver MUMPS.
\end{enumerate}
\end{algorithm}
\vspace{2mm}

\begin{remark}
The information of discontinuity or non-uniformity of the physical parameters has already been taken into account
in the devise of our preconditioner \eqref{sec4.4: eq_final_preconditioner}.
See the derivation of $\widehat{S}_{\phi}$, $M_p$, $L_p$ and $\widehat{\BF}$.
So we expect the performance to be robust with respect to the discontinuity of $\rho_h^n,\sigma_h^n$ and $\nu_h^n$.
Because $\gamma\left(\Div\BJ_h,\Div\boldsymbol{\varphi}_h\right)$ and
$\beta (\Div\Bu_h,\Div\Bv_h)$ are added in the discrete current density equation
and momentum equation, we call $\bbP$ the augmented
block preconditioner.
\end{remark}
\vspace{2mm}

\begin{remark}
Because $\sigma_h^n$ may have strong discontinuity, the condition of $\widetilde{\BQ}_\BJ$ will be
very bad. For the present we have not realized the efficient solver for $\widetilde{\BQ}_\BJ$.
So we use direct solver instead to mainly assess the effect of the whole block preconditioner $\bbP$.
Currently we are investigating and implementing the extension of the algorithm of auxiliary space preconditioning in
\cite{hip2007} to $\widetilde{\BQ}_\BJ$ on moving curved tetrahedral meshes in PHG.
\end{remark}

%%%%%%%%%%%%%%%%%%%%%%%%%%%%%%%%%%%%%%%%%%%
%%%%%%%%%%%%%%%%%%%%%%%%%%%%%%%%%%%%%%%%%%%
\section{Numerical experiments}\label{sec:numerics}
In this section, we present extensive tests to verify the desired properties of the proposed scheme
and preconditioner. In subsection~\ref{sec5.1}, we verified the computational accuracy and charge conservation
of the proposed high-order ALE finite element algorithm. In subsection~\ref{sec:5.2},
we investigate the stability and robustness of the proposed preconditioner
through the problem of drop of heavy fluid in the light fluid.
The finite element method and the solver are both implemented using the finite element software--Parallel Hierarchical Grid (PHG) \cite{zha22}. The 3D figures are generated by ParaView and right-preconditioned FGMRES without restart is used.

%%%%%%%%%%%%%%%%%%%%%%%%%%%%%%%%%%%%%%%%%%%%%%%%%%%
\subsection{Convergence and accuracy test}\label{sec5.1}
To verify the numerical accuracy of the proposed method, we consider the initial computational domain $\Omega_0=[0,1]^3$ and prescribe the imposed magnetic field as $\boldsymbol{B}=(0,1,0)^\top$. The corresponding manufactured solutions for the current density, velocity, electric potential and pressure are given by
\begin{align*}
\BJ =\,& \big(\sin (z+y),\, \cos (\pi z),\, \sin(x+y)\big)^\top,\\
\phi = \,& \sin(z),\\
\Bu =\,& \big(\sin (\pi x) \cos (\pi y),\, -\cos (\pi x)\sin (\pi y),\, 0\big)^\top,\\
p =\,& \frac{1}{4}\big(\cos (2\pi x) + \cos (2\pi y) \big).
\end{align*}
For simplicity, all parameters are set to be constants:
\[\sigma=Re=\alpha=\rho=\nu=1, \quad   \beta=0.5, \quad \gamma = 1.\]
Time-step length $\Delta t = 10^{-4}$ and terminal time $T_f=10^{-2}$ such that $N=100$. The initial guess of FGMRES iteration in every time step is set by the solution at previous time step and relative tolerance is set by $10^{-8}$.
For spatial discretization, velocity--pressure pair is approximated by the parametric $\boldsymbol{P}_{k+1}$--$P_k$ Taylor--Hood elements. The current density is discretized by the parametric $\mathbf{BDM}_k$ element, while the electric potential is approximated by the
parametric discontinuous Galerkin element of degree $k-1$ .
To verify the high-order accuracy of the proposed fully discrete
scheme \eqref{eq:fully_discrete_scheme}, we consider two cases with $k=1$ and $k=2$, respectively.
The corresponding convergence results are
presented in Table~\ref{tab:convergence_2order} and Table~\ref{tab:convergence_3order}.

\begin{figure}[!htbp]
	\centering
	\includegraphics[width=10cm]{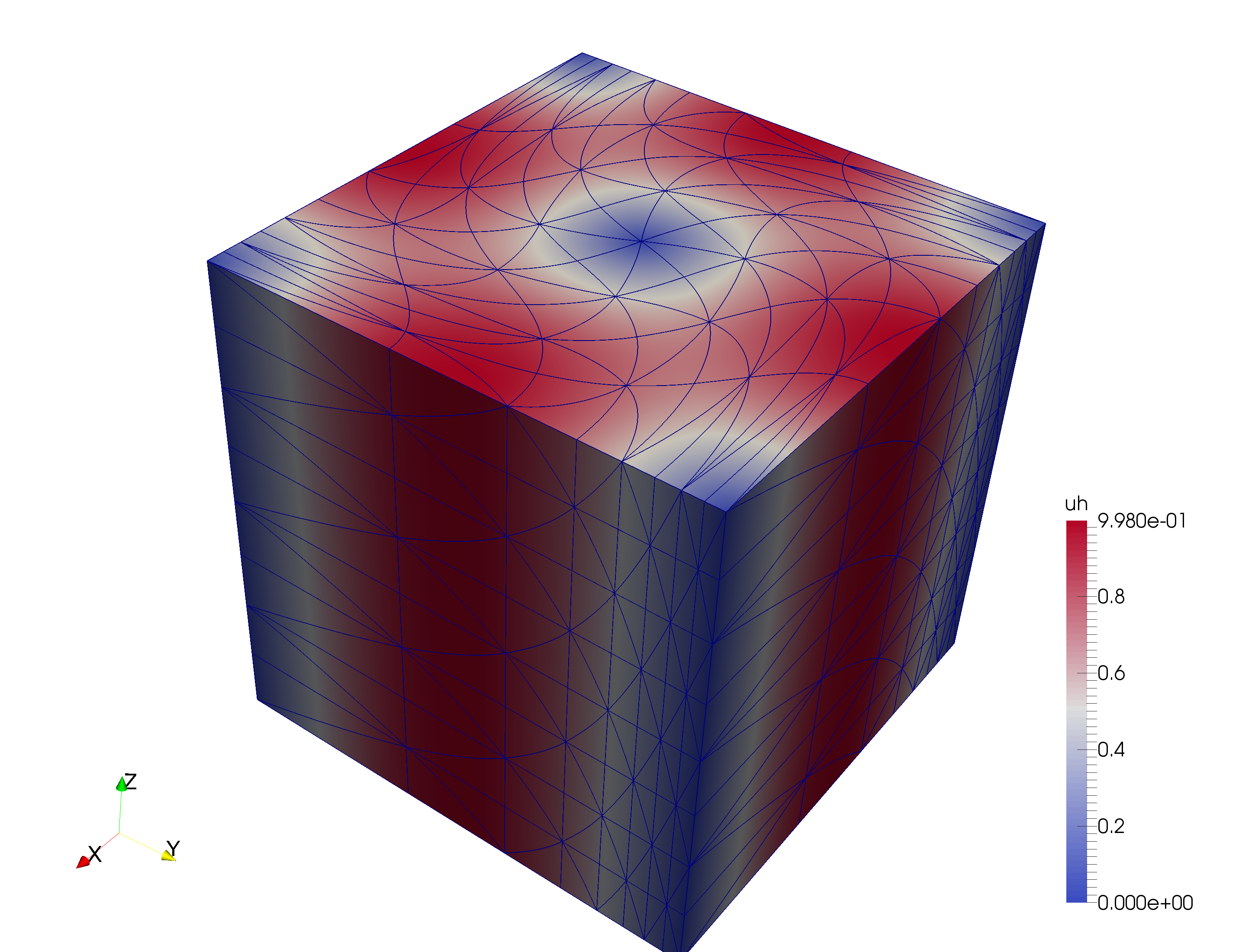}
	\caption{Curved mesh and 3D distribution of $|\Bu_h|$ on time $T_f=0.4$ with $k=1$.}
    \label{ex1:fig1}
\end{figure}

For the two finite element spaces with $k=1$ and $k=2$, the computed convergence rates agree well with the theoretical predictions.
For $k=1$, the proposed fully discrete scheme achieves the following second-order spatial convergence rates:
\[
   \left\|\boldsymbol{u}-\boldsymbol{u}_h^N\right\|_{1,\Omega_h^N}\approx O(h^2),
   \qquad
   \left\|p-p_h^N\right\|_{0,\Omega_h^N}\approx O(h^2),
   \qquad
   \left\|\boldsymbol{J}-\boldsymbol{J}_h^N\right\|_{\operatorname{div},\Omega_h^N}\approx O(h^2).
\]
When $k=2$, the convergence rates are further improved to three-order spatial convergence rates
\[
   \left\|\boldsymbol{u}-\boldsymbol{u}_h^N\right\|_{1,\Omega_h^N}\approx O(h^3),
   \qquad
   \left\|p-p_h^N\right\|_{0,\Omega_h^N}\approx O(h^3),
   \qquad
   \left\|\boldsymbol{J}-\boldsymbol{J}_h^N\right\|_{\operatorname{div},\Omega_h^N}\approx O(h^3).
\]
Moreover the last columns of Tables~\ref{tab:convergence_2order} and \ref{tab:convergence_3order} show that $\left\|\nabla \cdot \boldsymbol{J}_h^N\right\|_{0, \Omega_h^N} \approx 0$,
thereby numerically verifying the discrete charge conservation property.
Finally letting $T_f = 0.4, \Delta t = 10^{-3}$ and $k=1$, we show the curved tetrahedral mesh and distribution
of $|\Bu_h|$ in Fig.~\ref{ex1:fig1}.

\begin{table}[!htbp]
  \centering
  \fontsize{9pt}{12pt}\selectfont
  \caption{Convergence orders of the finite element solutions with $k=1$. (Section~\ref{sec5.1}).}
  \vspace{0.1in}
  \label{tab:convergence_2order}
  \begin{tabular}{@{\hspace{0.5em}}c@{\hspace{0.5em}}
                c@{\hspace{0.5em}}
                c@{\hspace{0.5em}}
                c@{\hspace{0.5em}}
                c@{\hspace{0.5em}}
                c@{\hspace{0.5em}}
                c@{\hspace{0.5em}}
                c@{\hspace{0.5em}}}
    \hline
    $h$ &
    $\left\|\boldsymbol{u}-\boldsymbol{u}_h^N\right\|_{1, \Omega_h^N}$ &
    order &
    $\left\|p-p_h^N\right\|_{0, \Omega_h^N}$ &
    order &
    $\left\|\boldsymbol{J}-\boldsymbol{J}_h^N\right\|_{\operatorname{div}, \Omega_h^N}$ &
    order &$\left\|\nabla \cdot \boldsymbol{J}_h^N\right\|_{0, \Omega_h^N}$  \\
    \hline
    0.8660 &3.40e-01 & -      &9.87e-02 & -      &4.80e-02 & -     &3.31e-10 \\
    0.4330 &1.82e-01 &0.91    &1.17e-01 &-0.24   &1.28e-02 &1.91   &5.64e-11 \\
    0.2165 &4.60e-02 &1.98    &2.44e-02 &2.26    &3.13e-03 &2.03   &2.96e-11 \\
    0.1083 &1.15e-02 &2.00    &4.76e-03 &2.36    &7.76e-04 &2.01   &8.34e-10 \\
    \hline
  \end{tabular}
\end{table}

\begin{table}[!htbp]
  \centering
  \fontsize{9pt}{12pt}\selectfont
  \caption{Convergence orders of the finite element solutions with $k=2$ (Section~\ref{sec5.1}).}
  \vspace{0.1in}
  \label{tab:convergence_3order}
  \begin{tabular}{@{\hspace{0.5em}}c@{\hspace{0.5em}}
                c@{\hspace{0.5em}}
                c@{\hspace{0.5em}}
                c@{\hspace{0.5em}}
                c@{\hspace{0.5em}}
                c@{\hspace{0.5em}}
                c@{\hspace{0.5em}}
                c@{\hspace{0.5em}}}
    \hline
    $h$ &
    $\left\|\boldsymbol{u}-\boldsymbol{u}_h^N\right\|_{1, \Omega_h^N}$ &
    order &
    $\left\|p-p_h^N\right\|_{0, \Omega_h^N}$ &
    order &
    $\left\|\boldsymbol{J}-\boldsymbol{J}_h^N\right\|_{\operatorname{div}, \Omega_h^N}$ &
    order &$\left\|\nabla \cdot \boldsymbol{J}_h^N\right\|_{0, \Omega_h^N}$
    \\
    \hline
    0.8660 &1.97e-01  & -      &1.31e-01  & -    &8.41e-03  & -     &1.71e-10 \\
    0.4330 &1.86e-02  &3.41    &8.60e-03  &3.93  &7.77e-04  &3.44   &1.67e-10 \\
    0.2165 &2.27e-03  &3.03    &8.78e-04  &3.29  &9.14e-05  &3.09   &2.53e-10 \\
    0.1083 &2.80e-04  &3.02    &9.80e-05  &3.16  &1.13e-05  &3.01   &2.69e-10 \\
    \hline
  \end{tabular}
\end{table}

%%%%%%%%%%%%%%%%%%%%%%%%%%%%%%%%%%%%%%%%%%%%%%%%%%%
\subsection{Performance test of the augmented block preconditioner}\label{sec:5.2}
\begin{figure}[!htbp]
	\centering
	\includegraphics[width=10cm]{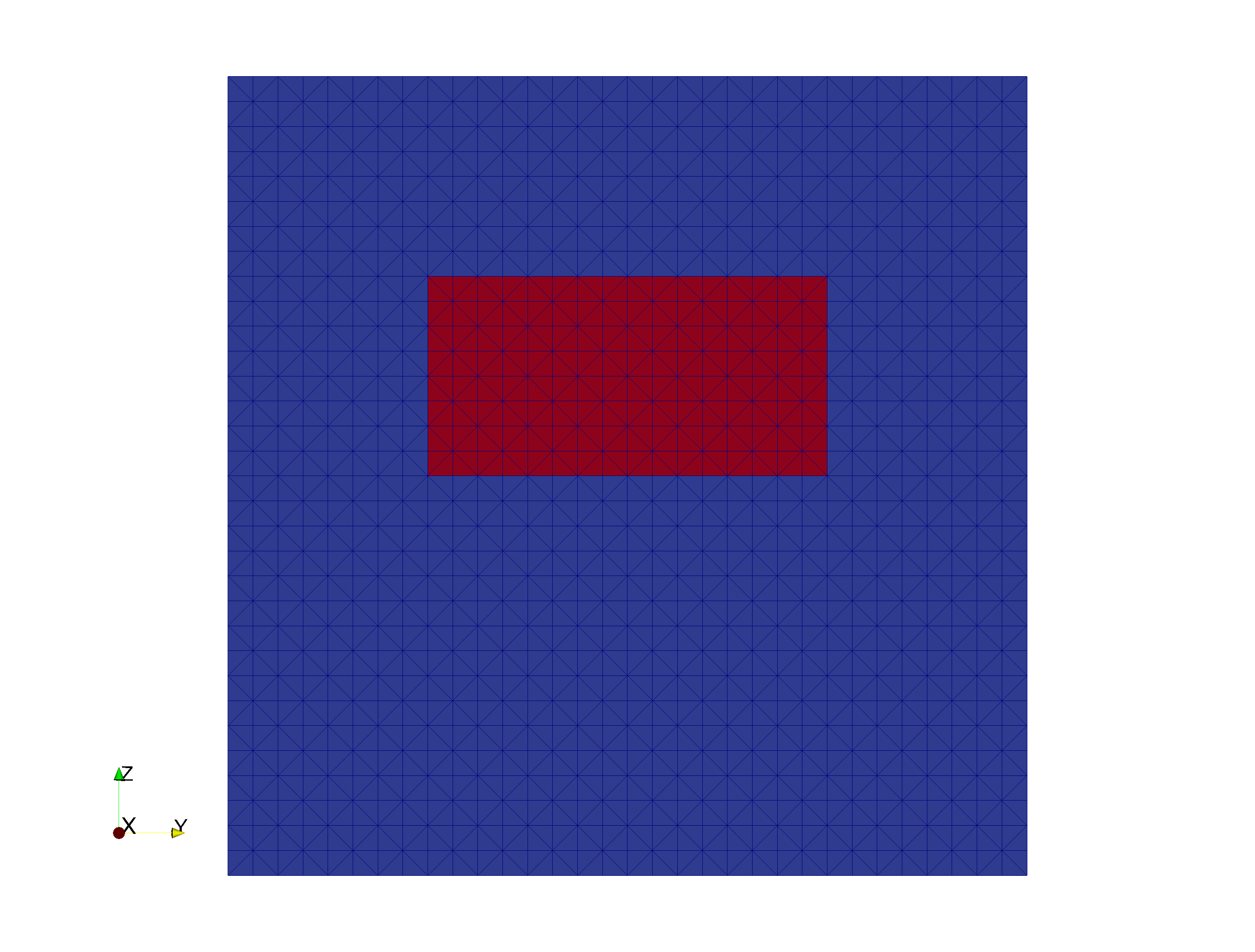}
	\caption{Initial distribution of the two phase on cross section $x = 1$.
    The red one represents heavy fluid (Section~\ref{sec:5.2}).}
    \label{ex2:fig1}
\end{figure}

Unlike single-phase MHD problems,
two-phase MHD problem considered here involves jumps in the physical parameters across the interface,
which increases the difficulty of devise of robust preconditioning.
To fully assess the robustness of the proposed preconditioner \eqref{sec4.4: eq_final_preconditioner}
with respect to main physical parameters and mesh refinement,
we use a 3D drop problem of heavy fluid into light fluid as benchmark.
In this section, we firstly summarize the fixed characteristic quantities in Table~\ref{tab:characteristic_scales}.
\begin{table}[H]
  \centering
  \caption{Characteristic quantities which are fixed for nondimensionalization (Section~\ref{sec:5.2}).}
  \label{tab:characteristic_scales}
  \begin{tabular}{lccc}
  \hline
  Characteristic quantity & Symbol & Value & Unit \\
  \hline
  Density
  & $\rho_c$ & $1000$ & $\mathrm{kg/m^3}$ \\
  Dynamic viscosity
  & $\nu_c$ & $10^{-3}$ & $\mathrm{Kg/m\cdot s}$ \\
  Electrical conductivity
  & $\sigma_c$ & $10^{4}$ & $\mathrm{S/m}$ \\
  Velocity
  & $u_c$ & $1$ & $\mathrm{m/s}$ \\
  \hline
\end{tabular}
\end{table}
The applied magnetic field is $(0,B_c,0)^\top$ and gravity $(0,0,-9.8\,\mathrm{m/s^2})^\top$.
Dimensionless initial domain is always set by
\[\Omega_{0} =  [0,2] \times[0,2]\times[-1,1],\]
The initial heavy fluid occupies the domain $\Omega_+(0) = [0.5,1.5]\times[0.5,1.5]\times[0,0.5]$,
which is marked by red in Fig.~\ref{ex2:fig1}.

Then initial material interface $\Gamma_0$ is set by $\partial\Omega_+(0)$.
Initial condition for $\Bu$ is set by $\Bu_0 = \mathbf{0}$.
For $\Bx_0\in\Omega_0$, initial density, electrical conductivity and viscosity coefficient are prescribed as
\begin{equation}
	\rho_0=
	\begin{cases}
		\rho_{\max}, \quad\Bx_0 \in \Omega_+(0),\\
		\rho_{\min}, \quad\Bx_0 \in \Omega_{-}(0),
	\end{cases}
\end{equation}
\begin{equation}
\sigma_0=
	\begin{cases}
		\sigma_{\max}, \quad\Bx_0 \in \Omega_+(0),\\
		\sigma_{\min}, \quad\Bx_0 \in \Omega_{-}(0),
	\end{cases}
\end{equation}
\begin{equation}
\nu_0=
	\begin{cases}
		\nu_{\max}, \quad\Bx_0 \in \Omega_+(0),\\
		\nu_{\min}, \quad\Bx_0 \in \Omega_{-}(0),
	\end{cases}
\end{equation}
In the following test, for easy presentation we set $\rho_{\max}=\sigma_{\max}=\nu_{\max}=1$
and vary $\rho_{\min}, \sigma_{\min}, \nu_{\min}$ to test the performance of the block preconditioner.
No-slip boundary condition for velocity and insulating boundary condition for current density are used
\begin{equation*}
\boldsymbol{u} = \boldsymbol{0}, \quad \BJ\cdot \Bn = 0 \qquad \text{on}~~\partial\Omega.
\end{equation*}

The initial guess of FGMRES iteration in every time step is set by solutions at previous time step and the relative tolerance is $10^{-6}$.
We use two time-step length $\Delta t = 0.1$ and $0.01$
to show the preconditioning effect for different $\Delta t$.
\textit{The average FGMRES iteration numbers will be reported,
namely the total iteration numbers divided by the total time steps
(we round off the value to remove its fractional part)}.
In this section, $k=1$ is used for convenience.
Table~\ref{table_Grid_information} presents the numbers of degrees of freedom for current density $\BJ_h$,
velocity $\Bu_h$, electric potential $\phi$ and pressure $p$ at different mesh sizes.

\begin{table}[!htbp]
  \centering
  \caption{Mesh information and the numbers of degree of freedoms with $k=1$ (Section~\ref{sec:5.2}).}
  \label{table_Grid_information}
  \begin{tabular}{|c|c|c|c|c|c|}\hline
  Mesh &$h$   &DOFs of $\BJ_h$  &DOFs of $\Bu_h$   &DOFs of $\phi_h$  &DOFs of $p_h$  \\\hline
  M1  &0.8660 &2592              &2187             &384                &125      \\
  M2  &0.4330 &19584             &14739            &3072               &729      \\
  M3  &0.2165 &152064            &107811           &24576              &4913     \\
  M4  &0.1083 &1198080           &823875           &196608             &35937    \\\hline
  \end{tabular}
\end{table}
\subsubsection{Sensitivity of the preconditioner for augmentation parameter $\beta$ and $\gamma$}\label{secpresub1}
In this subsection, we assess the robustness of the proposed solver with respect to the augmentation parameters
$\beta$ and $\gamma$. We set $L_c = 0.1\,\mathrm{m}$ and $B_c =1\,\mathrm{T}$ such that
\[Re = 10^5,~~\alpha = 1,~~\BB=(0,1,0)^\top,~~\Bf_u = (0,0,-0.98)^\top.\]
Other physical parameters are fixed by
\[\rho_{\min}=0.5,\quad \sigma_{\min}=10^{-3},\quad \nu_{\min}=0.25.\]
Terminal time $T_f=0.2$. Using different $\gamma$,
we report the average numbers of FGMRES iterations on
meshes M1--M3 for $\beta=0.1, 0.25, 0.5, 0.75$ and $1$.
The results are presented in Tables~\ref{sec521}--\ref{sec521b}.

\begin{table}[!htbp]
  \centering
  \caption{Set $\gamma = 0.1$: average FGMRES iteration numbers for different $\beta$ and $\Delta t$ (Subsection~\ref{secpresub1}).}
  \label{sec521}
  \begin{tabular}{|c|c|c|c|c|c|c|c|c|c|c|}\hline
  $\Delta t$
  &\multicolumn{5}{c|}{$0.1$ (2 steps)}
  &\multicolumn{5}{c|}{$0.01$ (20 steps)}\\\hline
  \diagbox[width=2cm]{Mesh}{$\beta$}
     &0.1 &0.25 &0.5 &0.75 &1.0  &0.1 &0.25 &0.5 &0.75 &1.0                       \\\hline
  M1 &13  &13 &14  &14  &14  &9  &9  &9 &9 &9  \\
  M2 &12  &11 &10  &9   &9   &9  &9  &9 &9 &9  \\
  M3 &13  &11 &10  &9   &9   &10 &9  &9 &9 &9  \\\hline
  \end{tabular}
\end{table}

\begin{table}[!htbp]
  \centering
  \caption{Set $\gamma = 0.5$: average FGMRES iteration numbers for different $\beta$ and $\Delta t$ (Subsection~\ref{secpresub1}).}
  \label{sec521a}
  \begin{tabular}{|c|c|c|c|c|c|c|c|c|c|c|}\hline
  $\Delta t$
  &\multicolumn{5}{c|}{$0.1$ (2 steps)}
  &\multicolumn{5}{c|}{$0.01$ (20 steps)}\\\hline
  \diagbox[width=2cm]{Mesh}{$\beta$}
     &0.1 &0.25 &0.5 &0.75 &1.0  &0.1 &0.25 &0.5 &0.75 &1.0                       \\\hline
  M1 &13  &13 &14  &14  &14  &9 &9  &9 &9 &8  \\
  M2 &12  &10 &9   &9   &8   &9 &9  &9 &8 &8  \\
  M3 &12  &10 &9   &8   &8   &9 &9  &8 &8 &8  \\\hline
  \end{tabular}
\end{table}

\begin{table}[!htbp]
  \centering
  \caption{Set $\gamma = 1.0$: average FGMRES iteration numbers for different $\beta$ and $\Delta t$ (Subsection~\ref{secpresub1}).}
  \label{sec521b}
  \begin{tabular}{|c|c|c|c|c|c|c|c|c|c|c|}\hline
  $\Delta t$
  &\multicolumn{5}{c|}{$0.1$ (2 steps)}
  &\multicolumn{5}{c|}{$0.01$ (20 steps)}\\\hline
  \diagbox[width=2cm]{Mesh}{$\beta$}
     &0.1 &0.25 &0.5 &0.75 &1.0  &0.1 &0.25 &0.5 &0.75 &1.0                       \\\hline
  M1 &12  &13 &14  &14  &14  &9 &9  &9 &8 &8  \\
  M2 &12  &10 &9   &8   &8   &9 &9  &8 &8 &8  \\
  M3 &12  &10 &9   &8   &8   &9 &9  &9 &8 &8  \\\hline
  \end{tabular}
\end{table}

Overall, for given $\beta$ and $\gamma$ the average FGMRES iteration counts are nearly independent of mesh refinement.
Given $\gamma >0$, larger values of $\beta$ yield slightly better results.
For example, when $\Delta t = 0.1$, they generally decrease as $\beta$ increases from $0.1$ to $1$.
In addition, smaller $\Delta t$ results in less iteration counts.
Moreover, for fixed $\Delta t$ and $\beta$, slightly smaller iteration count is obtained with larger $\gamma$.
These results demonstrate the robustness of the proposed preconditioner with
respect to the mesh size for moderate $\beta$ and $\gamma$.
Based on the data in Table~\ref{sec521}--\ref{sec521b},
in the following test, we will set $\beta =0.5$, $\gamma = 1$
and mainly use $\Delta t = 1 \times 10^{-2}$ to verify the robustness of the preconditioner.

\subsubsection{Sensitivity of the solver for $\sigma_{\min}$}\label{secpresub3}
\begin{figure}[!htbp]
	\centering
	\includegraphics[width=9cm]{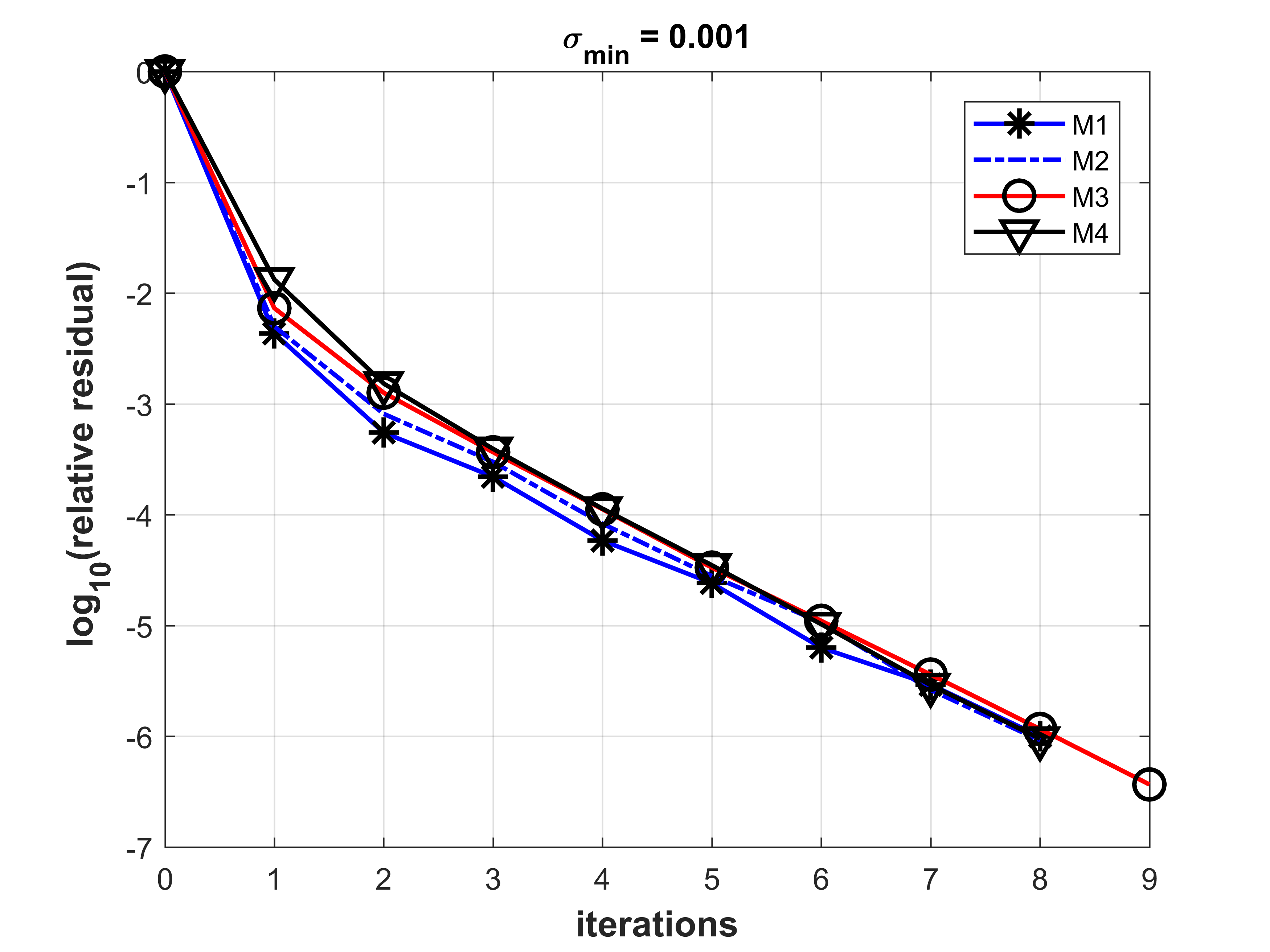}
	\caption{Decreasing trend of residual of FGMRES for $\sigma_{\min} = 10^{-3}$ (Subsection~\ref{secpresub3}).}
    \label{ex2:fig2}
\end{figure}

In this subsection, we assess the robustness of the proposed preconditioner with respect to the parameter $\sigma_{\min}$. We set $L_c = 0.1\,\mathrm{m}$ and $B_c =1\,\mathrm{T}$ such that
\[Re = 10^5,~~\alpha = 1,~~\BB=(0,1,0)^\top,~~\Bf_u = (0,0,-0.98)^\top.\]
We fix other parameters by
\[\rho_{\min}=0.5,\quad \nu_{\min}=0.25, \quad \beta =0.5, \quad\gamma = 1.\]
The terminal time is $T_f=0.2$. As shown in Table~\ref{sec5.2.3: sigma_minnew}, the average FGMRES iteration numbers are robust with respect to both the mesh size and $\sigma_{\min}$. The average iteration numbers remain approximately 8 on meshes M1--M4 for different $\sigma_{\min}$.
Using $\sigma_{\min} = 10^{-3}$, we plot the decreasing trend of FGMRES in Fig.~\ref{ex2:fig2}, which shows the mesh independence of the proposed block preconditioner.

\begin{table}[!htbp]
  \centering
  \caption{Average FGMRES iteration numbers for different $\sigma_{\min}$
  with $\Delta t = 1 \times 10^{-2}$ (Subsection~\ref{secpresub3}).}
  \label{sec5.2.3: sigma_minnew}
   \begin{tabular}{|c|c|c|c|c|c|}\hline
  \diagbox[width=2cm]{Mesh}{$\sigma_{\min}$}
      &0.1  & $10^{-2}$ &$10^{-3}$  &$10^{-4}$ &$10^{-6}$       \\\hline
  M1  &8    &8    &9     &8   &8      \\
  M2  &8    &8    &8     &8   &8                      \\
  M3  &8    &8    &8     &8   &8                       \\
  M4  &8    &8    &8     &9   &8                      \\\hline
  \end{tabular}
\end{table}

\subsubsection{Sensitivity of the solver for $\rho_{\min}$}\label{secpresub2}
\begin{figure}[!htbp]
	\centering
	\includegraphics[width=9cm]{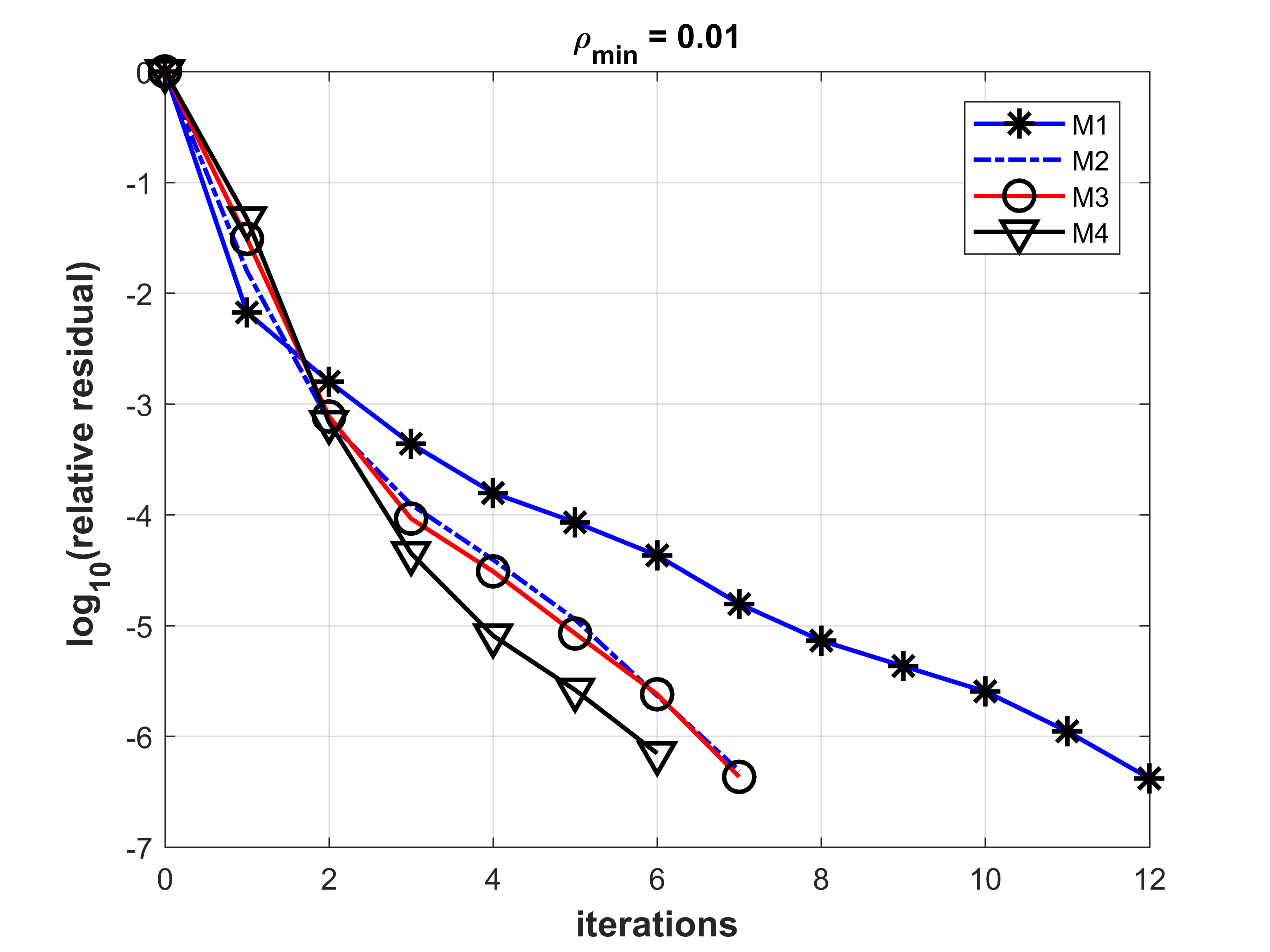}
	\caption{Decreasing trend of residual of FGMRES for $\rho_{\min} = 0.01$ (Subsection~\ref{secpresub2}).}
    \label{ex2:fig3}
\end{figure}

In this subsection we continue to
investigate the effect of $\rho_{\min}$ on the average number of FGMRES iterations for different mesh sizes.
We set $L_c = 0.1\,\mathrm{m}$ and $B_c =1\,\mathrm{T}$ such that
\[Re = 10^5,~~\alpha = 1,~~\BB=(0,1,0)^\top,~~\Bf_u = (0,0,-0.98)^\top.\]
We fix other physical parameters by
\[\sigma_{\min}=10^{-3},\quad \nu_{\min}=0.25, \quad \beta =0.5, \quad\gamma = 1.\]
Terminal time is $T_f=0.2$.
\begin{table}[!htbp]
  \centering
  \caption{Average FGMRES iteration numbers for different
  $\rho_{\min}$ with $\Delta t = 1 \times 10^{-2}$ (Subsection~\ref{secpresub2}).}
  \label{sec5.2.2: rho_min}
 \begin{tabular}{|c|c|c|c|c|}\hline
  \diagbox[width=2cm]{Mesh}{$\rho_{\min}$}
      &0.5 &0.1  &0.01    &0.001         \\\hline
  M1  &9  &10 &11     &11            \\
  M2  &8  &8  &7      &6            \\
  M3  &8  &8  &7      &6            \\
  M4  &8  &8  &6      &6            \\\hline
  \end{tabular}
\end{table}
As shown in Table~\ref{sec5.2.2: rho_min}, the average FGMRES iteration numbers are robust with respect to both the mesh size and $\rho_{\min}$. Specifically, the average iteration numbers remain approximately 10 on meshes M1--M4.
The average iteration numbers are stable as $\rho_{\min}$ decreases from $0.5$ to $0.001$. Also, Fig.~\ref{ex2:fig3}
shows the decreasing trend of residual of FGMRES iteration.

\subsubsection{Sensitivity of the solver for $\nu_{\min}$}\label{secpresub4}
In this subsection, we assess the robustness of the proposed preconditioner with respect
to the parameter $\nu_{\min}$. We set $L_c = 0.1\,\mathrm{m}$ and $B_c =1\,\mathrm{T}$ such that
\[Re = 10^5,~~\alpha = 1,~~\BB=(0,1,0)^\top,~~\Bf_u = (0,0,-0.98)^\top.\]
We fix other physical parameters by
\[\rho_{\min}=0.5,\quad \sigma_{\min}=10^{-3}, \quad \beta =0.5, \quad\gamma = 1.\]
The terminal time is $T_f=0.2$. We report the average FGMRES iterations numbers,
namely the total iteration numbers divided by the total time steps.

\begin{table}[!htbp]
  \centering
  \caption{Average FGMRES iteration numbers for different $\nu_{\min}$
  with $\Delta t = 1 \times 10^{-2}$ (Subsection~\ref{secpresub4}).}
  \label{sec5.2.4: nu_min}
  \begin{tabular}{|c|c|c|c|c|}\hline
  \diagbox[width=2cm]{Mesh}{$\nu_{\min}$}
      &0.25 &0.1 &0.01  &0.001  \\\hline
  M1  &9    &9  &9       &9             \\
  M2  &8    &8  &8       &8              \\
  M3  &8    &8  &8       &8              \\
  M4  &8    &8  &8       &8            \\\hline
  \end{tabular}
\end{table}

\begin{figure}[!htbp]
	\centering
	\includegraphics[width=9cm]{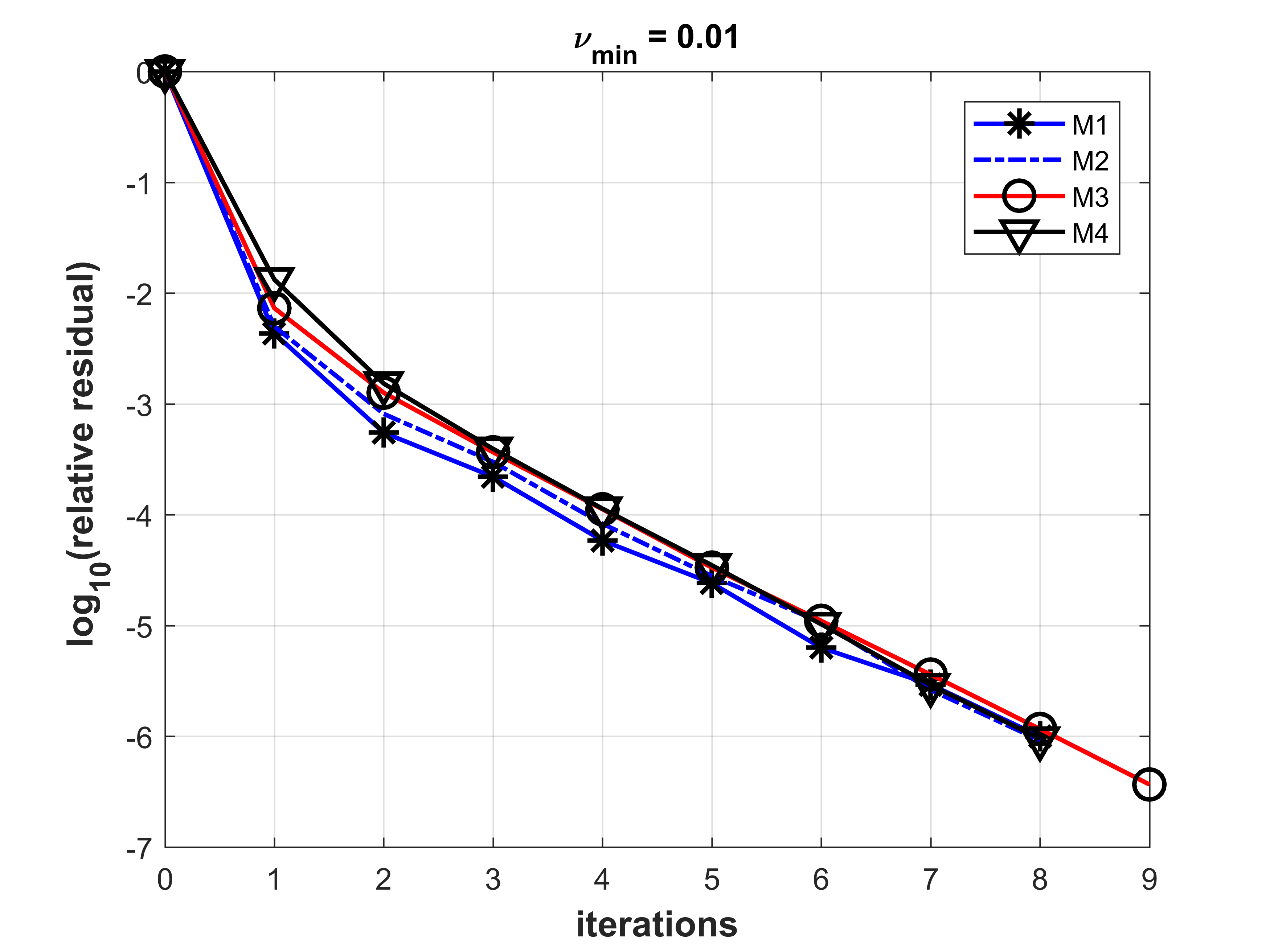}
	\caption{Decreasing trend of residual of FGMRES for $\nu_{\min} = 0.01$ (Subsection~\ref{secpresub4}).}
    \label{ex2:fig4}
\end{figure}

As shown in Table~\ref{sec5.2.4: nu_min}, the average FGMRES iteration numbers are robust with respect to
both the mesh refinement and $\nu_{\min}$.
Moreover, the average iteration numbers remain nearly unchanged as $\nu_{\min}$ decreases from $0.25$ to $0.001$.
Fig.~\ref{ex2:fig4} shows the decreasing trend of residual of FGMRES for $\nu_{\min} = 0.01$.
Again on can observe stable reduction of the residual.

\subsubsection{Sensitivity of the solver for $Re$ and $\alpha$}\label{secpresub5}
In the preceding subsections, we have verified the robustness of the proposed preconditioner with respect to the discontinuity of physical parameters $\rho$, $\sigma$, and $\nu$. In this subsection, we assess the robustness of the proposed preconditioner with respect to the Reynolds number $Re$ and the coupling number $\alpha$.
We fix the following parameters
\[\rho_{\min}=0.5,\quad \sigma_{\min}=10^{-3},\quad \nu_{\min}=0.25, \quad \beta =0.5, \quad\gamma = 1.\]
Recall the definition of Reynolds number, coupling number and dimensionless body force
\[Re  =  \frac{\rho_cu_cL_c}{\nu_c},\quad
\alpha = \frac{\sigma_cB_c^{2}L_c}{\rho_cu_c}, \quad
\Bf_u = \frac{L_c}{u_c^2}(0,0,-9.8\,\mathrm{m/s^2})^\top.\]
The terminal time is $T_f = 0.2$. We consider six different physical cases
\begin{enumerate}
  \item [\textbf{C1.}] $L_c=0.01\,\mathrm{m}, B_c = 1\,\mathrm{T}$ such that
  \[Re = 10^4, \alpha = 0.1, \Bf_u = (0,0,-9.8\times 10^{-2}).\]

  \item [\textbf{C2.}] $L_c=0.1\,\mathrm{m}, B_c = 1\,\mathrm{T}$ such that
  \[Re = 10^5, \alpha = 1, \Bf_u = (0,0,-0.98).\]

  \item [\textbf{C3.}] $L_c=1\,\mathrm{m}, B_c = 1\,\mathrm{T}$ such that
  \[Re = 10^6, \alpha = 10, \Bf_u = (0,0,-9.8).\]

  \item [\textbf{C4.}] $L_c=0.01\,\mathrm{m}, B_c = 10\,\mathrm{T}$  such that
  \[Re = 10^4, \alpha = 10, \Bf_u = (0,0,-9.8\times 10^{-2}).\]

  \item [\textbf{C5.}] $L_c=0.1\,\mathrm{m}, B_c = 10\,\mathrm{T}$ such that
  \[Re = 10^5, \alpha = 100, \Bf_u = (0,0,-0.98).\]

  \item [\textbf{C6.}] $L_c=1\,\mathrm{m}, B_c = 10\,\mathrm{T}$ such that
  \[Re = 10^6, \alpha = 1000, \Bf_u = (0,0,-9.8).\]
\end{enumerate}
The corresponding results are presented in Table~\ref{sec5.2.5: small_t}-\ref{sec5.2.5: large_t}.

\begin{table}[!htbp]
  \centering
  \caption{Average FGMRES iteration numbers for different physical cases with $\Delta t = 0.1$ (Subsection~\ref{secpresub5}).}
  \label{sec5.2.5: small_t}
 \begin{tabular}{|c|c|c|c|c|c|c|}\hline
  Mesh &C1  &C2      &C3  &C4  &C5 &C6 \\\hline
  M1   &13  &14    &16    &15    &25   &73   \\
  M2   &9   &9     &12    &12    &25   &78   \\
  M3   &8   &9     &11    &11    &25   &74   \\
  M4   &8   &8     &11    &11    &23   &71   \\\hline
  \end{tabular}
\end{table}

\begin{table}[!htbp]
  \centering
  \caption{Average FGMRES iteration numbers for different physical cases with $\Delta t = 0.01$ (Subsection~\ref{secpresub5}).}
  \label{sec5.2.5: large_t}
  \begin{tabular}{|c|c|c|c|c|c|c|}\hline
  Mesh &C1  &C2  &C3  &C4  &C5 &C6 \\\hline
  M1   &8   &9   &10  &10  &13 &27   \\
  M2   &7   &8   &10  &10  &13 &28   \\
  M3   &7   &8   &10  &10  &13 &28   \\
  M4   &7   &8   &11  &10  &13 &26 \\\hline
  \end{tabular}
\end{table}

\begin{figure}[!htbp]
	\centering
	\includegraphics[width=9cm]{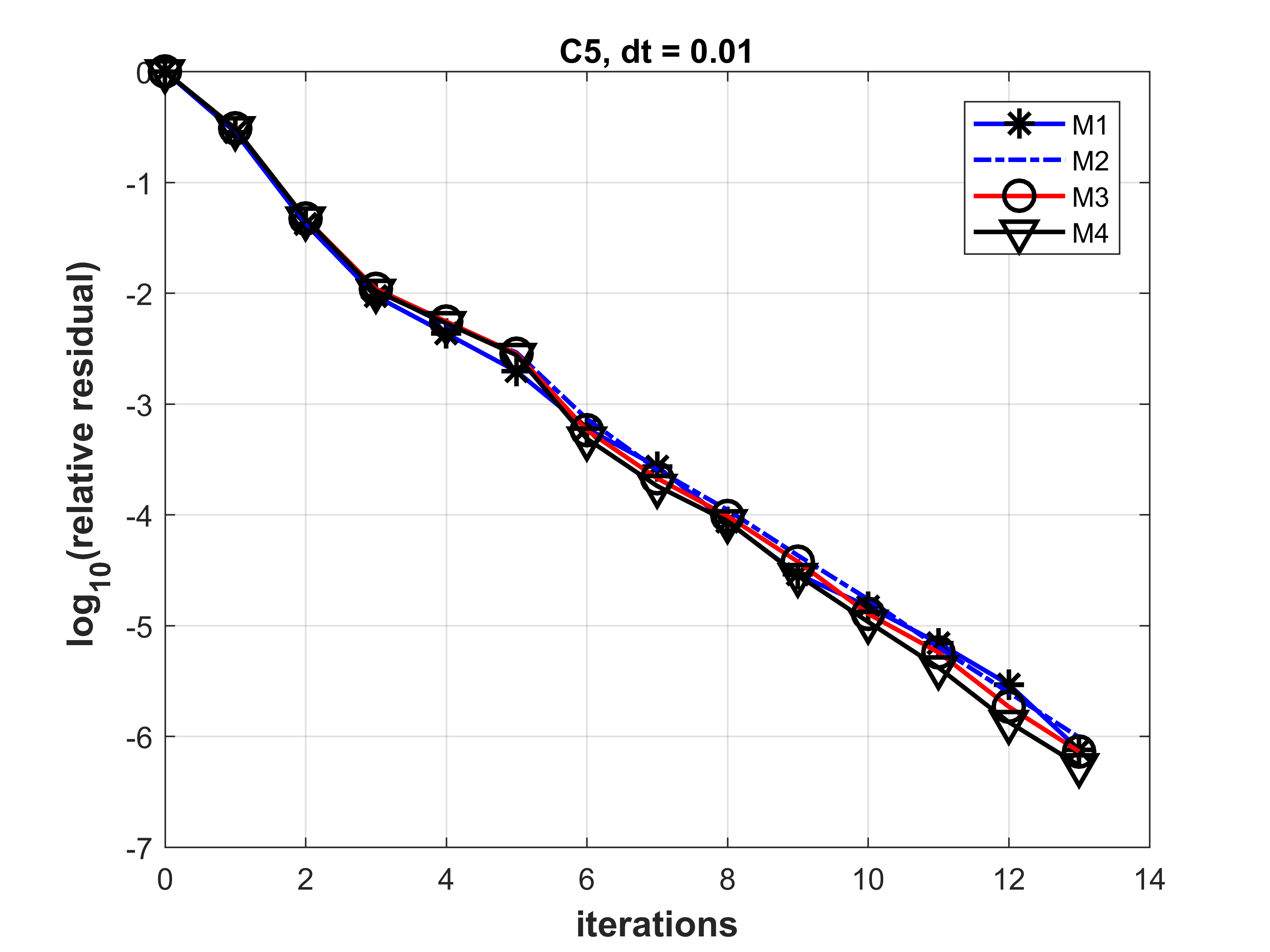}
	\caption{Decreasing trend of residual of FGMRES for C5 and $\Delta t = 0.01$ (Subsection~\ref{secpresub5}).}
    \label{ex2:fig5}
\end{figure}

Firstly, from both tables, we find the preconditioner is robust with respect to grid refinement.
For fixed mesh and physical case, smaller values of $\Delta t$ results in fewer iteration numbers.
However our preconditioner performance starts to deteriorate for very large coupling number and large time-step length.
For case C6, when $\Delta t = 0.1$, average iteration numbers exceed 70 steps.
But when $\Delta t = 0.01$, it reduces to roughly 28 steps.
In future work, we will investigate how to improve the preconditioning effect for large coupling numbers.
Finally, in Fig.~\ref{ex2:fig5}, we plot the decreasing trend of residual for case C5 with $\Delta t = 0.01$.

\subsubsection{Simulation of the drop process using the ALE solver}\label{secpresub6}
\begin{figure}[!htbp]
	\centering
	\includegraphics[width=10cm]{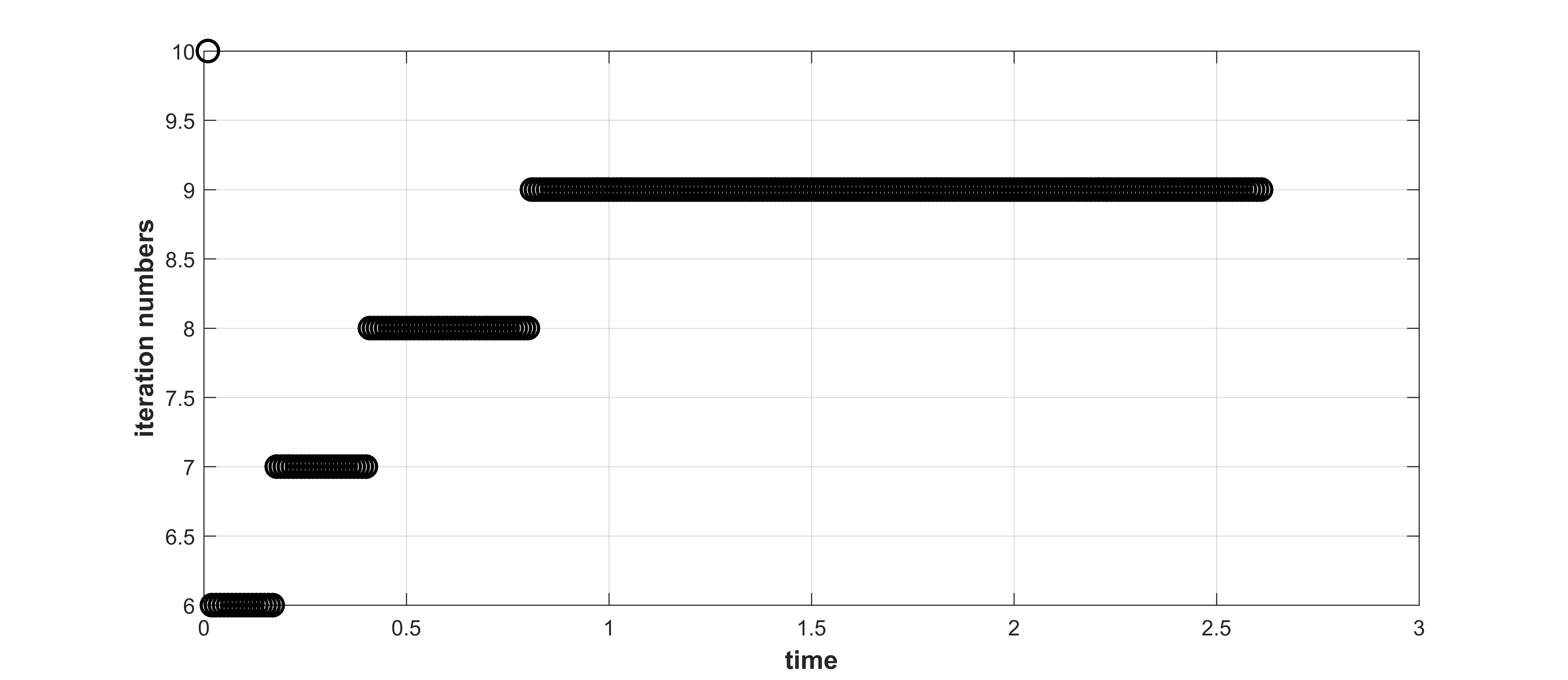}
	\caption{FGMRES iteration numbers in every time step (Subsection~\ref{secpresub6}).}
    \label{ex2:fig6}
\end{figure}

\begin{figure}[!htbp]
	\centering
	\includegraphics[width=6cm]{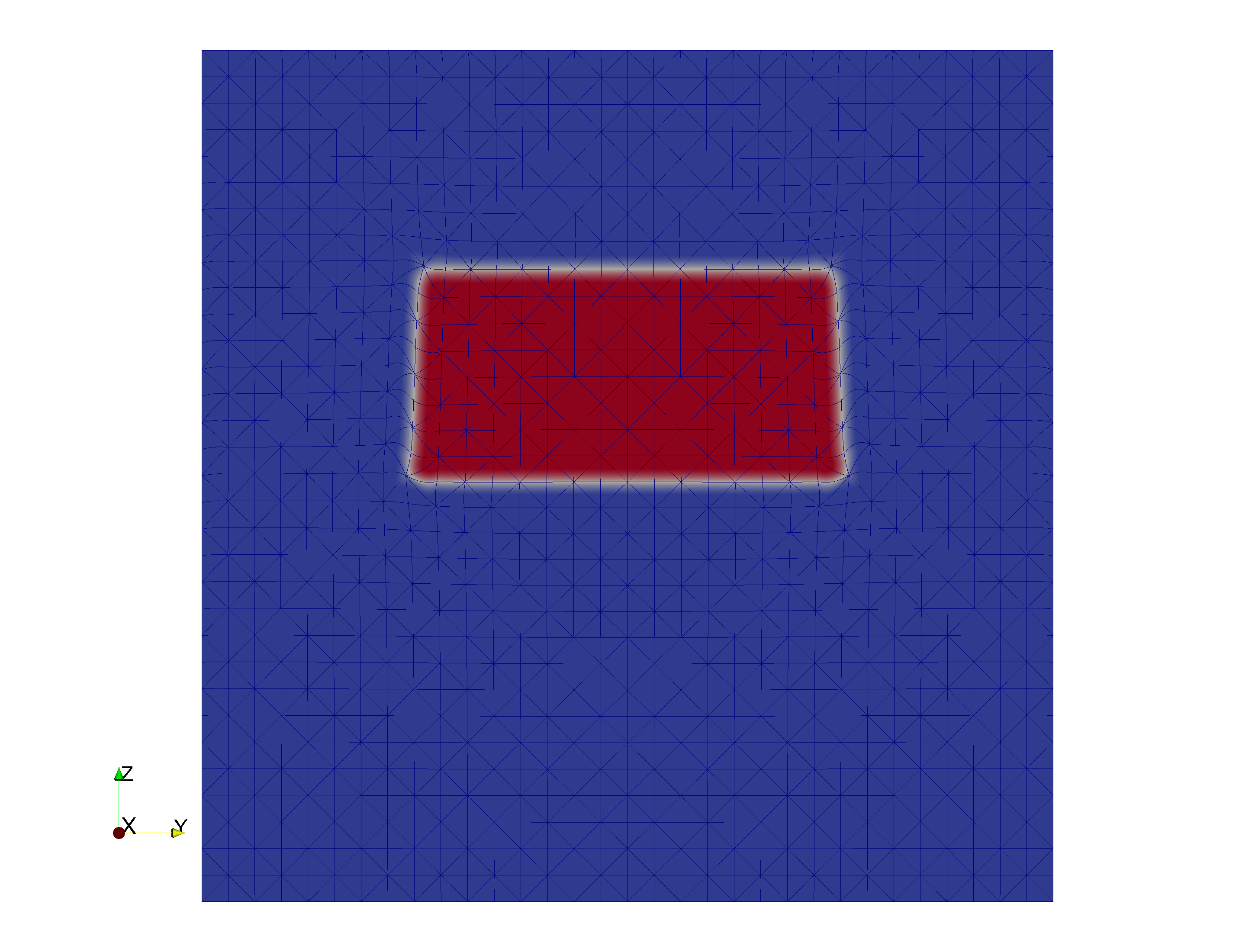}
    \includegraphics[width=6cm]{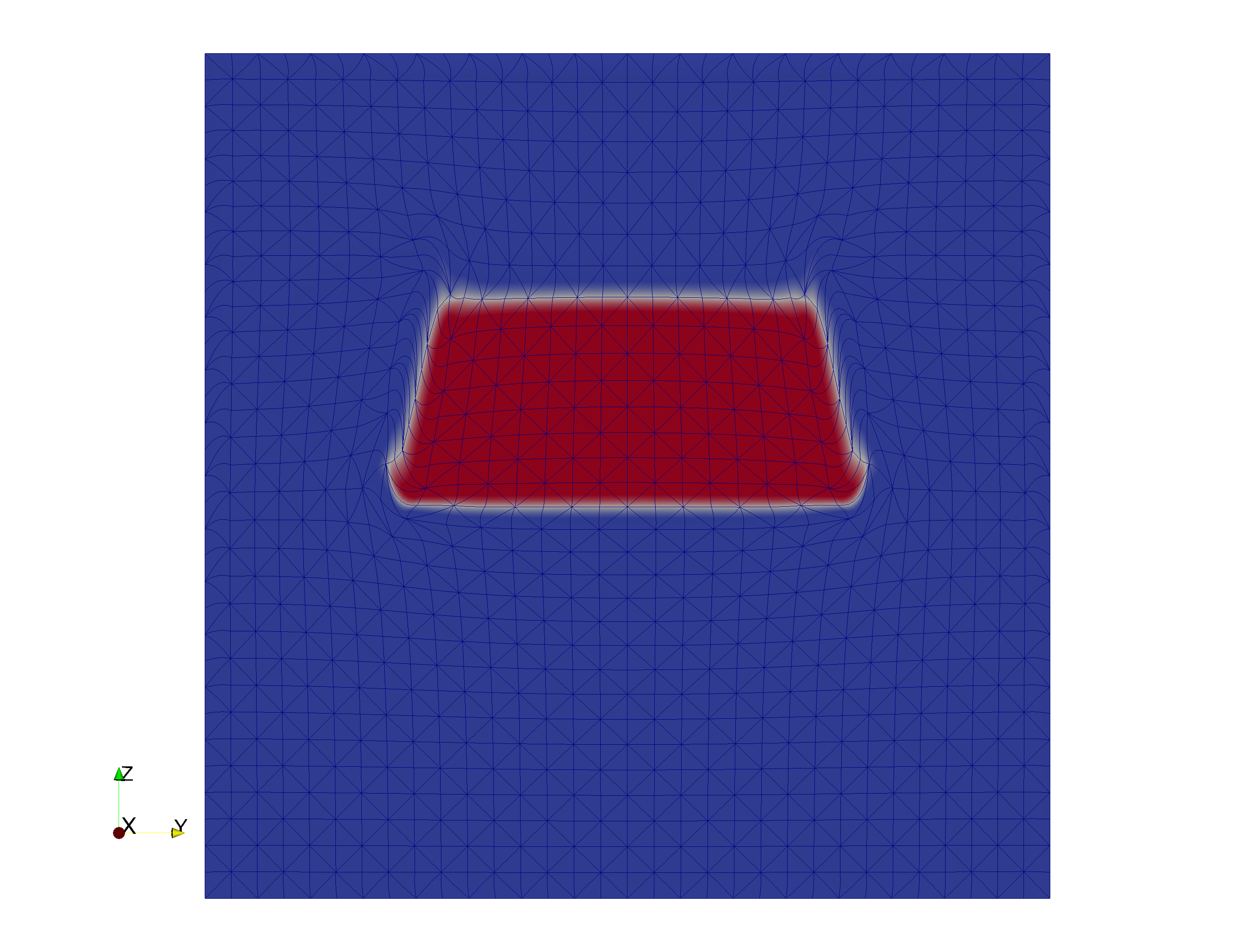}
	\caption{Material distribution on time $t=1.13$ and $t=2.6$ in the vicinity of cross section $x=0.5$
    (Subsection~\ref{secpresub6}).}
    \label{ex2:fig7}
\end{figure}

\begin{figure}[!htbp]
	\centering
	\includegraphics[width=6cm]{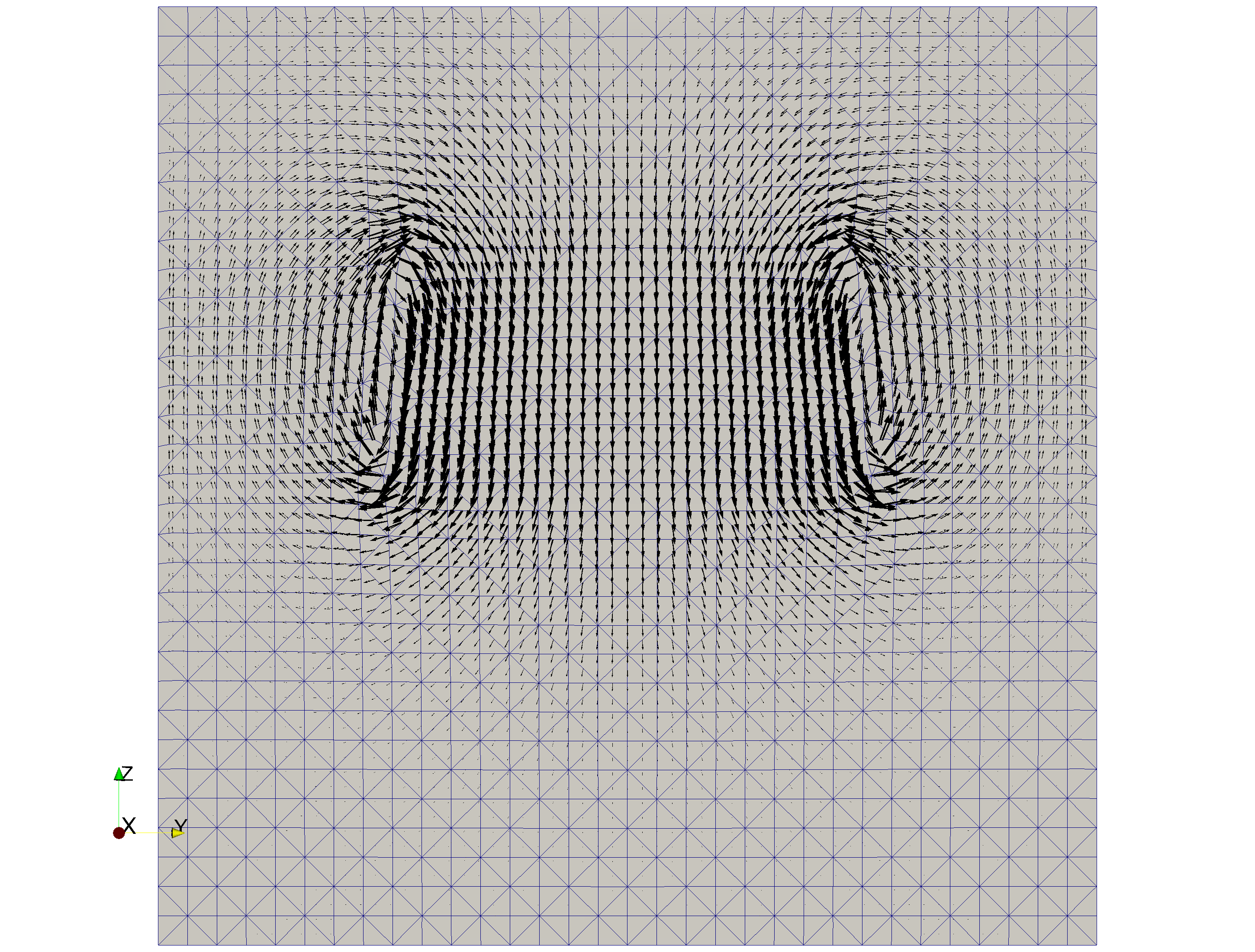}
    \includegraphics[width=6cm]{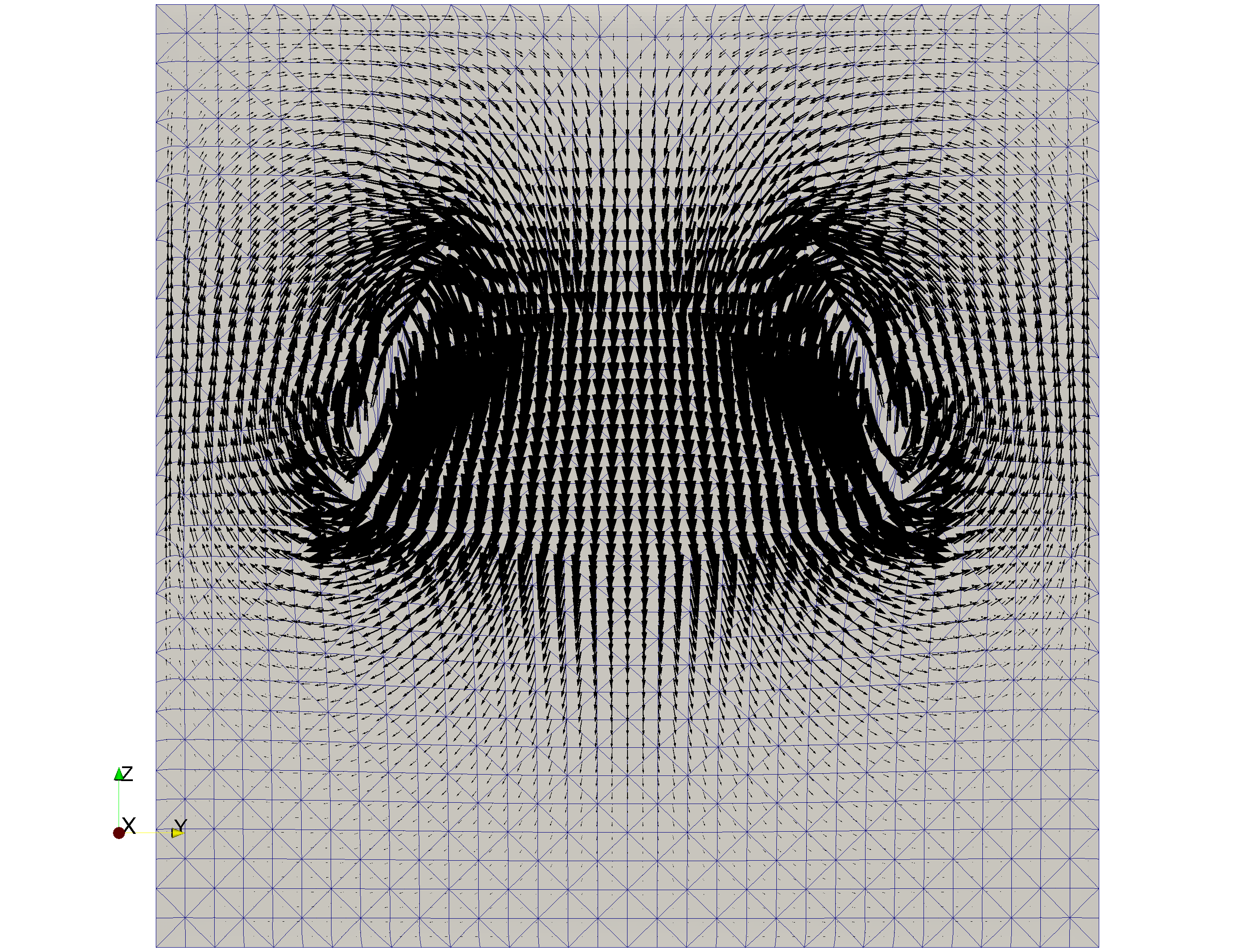}
	\caption{Velocity vector distribution
     on time $t=1.13$ and $t=2.6$ in the vicinity of cross section $x=0.5$
     (Subsection~\ref{secpresub6}).}
    \label{ex2:fig8}
\end{figure}

In this subsection, we use the mesh M4 and
set $L_c = 0.01\,\mathrm{m}$ and $B_c =1\,\mathrm{T}$ to simulate the drop process using
the proposed ALE finite element method.
In this case $Re = 10^4,\alpha = 0.1,\BB=(0,1,0)^\top$ and $\Bf_u = (0,0,-9.8\times 10^{-2})^\top$.
Other parameters are set by
\[\rho_{\min}=0.5,\quad \sigma_{\min}=0.1,\quad \nu_{\min}=0.25, \quad \beta =0.5, \quad\gamma = 1.\]
The terminal time $T_f=2.6$ and $\Delta t = 0.01$. Firstly, in Fig.~\ref{ex2:fig6}, we give the concrete
iteration numbers in every time step. Fig.~\ref{ex2:fig7} shows the distribution of the two phase on time $t=1.13$
and $t=2.6$ in the vicinity of cross section $x=0.5$. Fig.~\ref{ex2:fig8} gives the velocity vector distribution
on time $t=1.13$ and $t=2.6$ in the vicinity of cross section $x=0.5$.

\section{Conclusion}\label{sec:conclusion}
In this work, we mainly study efficient preconditioning for two-phase incompressible inductionless MHD equations.
We employ a high-order, charge-conservative, monolithic ALE mixed finite element method as base scheme.
All physical variables are solved simultaneously at each time step.
The resulting linear algebraic system is a large-scale double saddle-point system.
To solve it efficiently, we developed an augmented block preconditioner based on constraint preconditioning, weighted augmentation and approximate Schur complement strategy. Numerical experiments verify the spatial accuracy and charge conservation of the high-order ALE method. Extensive tests are given to demonstrate the robustness of the preconditioner with respect to mesh refinement and variations in the main physical parameters.
In future work, we will explore more efficient mesh smoothing strategy and
devise energy stable ALE method with surface tension.
\bibliographystyle{amsplain}

\end{document}